\documentclass[12pt]{amsart}
\usepackage{amsmath,amssymb,amscd,array, mathrsfs }

\usepackage{amsmath,amscd,amssymb,amsthm,array}
\usepackage{color}

\usepackage{amsmath,amssymb,mathrsfs,amsthm, tikz-cd,mathrsfs}
\usepackage{comment}
\def\url#1{\expandafter\s

\tring\csname #1\endcsname}

\def\mmat #1,#2,#3,#4,{\text{\small\arraycolsep=3pt $
\begin{pmatrix}#1&#2\\#3&#4\end{pmatrix}$}}

\usepackage{hyperref}
\usepackage{capt-of}

\usepackage{multirow}

\usepackage{comments}
\newComments\NG{Nico}{red}
\newComments\QEh{QEh}{blue}

\def\mmat #1,#2,#3,#4,{\text{\small\arraycolsep=3pt $
\begin{pmatrix}#1&#2\\#3&#4\end{pmatrix}$}}

\usepackage{lscape}
\usepackage{tikz-cd}
\usepackage{enumerate}

\usepackage{DLdef1}
\usepackage{multicol}

\def\mmat #1,#2,#3,#4,{\text{\small\arraycolsep=3pt $
\begin{pmatrix}#1&#2\\#3&#4\end{pmatrix}$}}

\renewcommand {\ssbegin}[2][*]
 {\refstepcounter{subsection}%
\if#1*
\addcontentsline{toc}{subsection}{\thesubsection.\hskip 1pc #2}%
\else
\addcontentsline{toc}{subsection}{\thesubsection.\hskip 1pc #2. #1}%
\fi
 \def \secno {\gdef \secno {}{\ssecfont
\thesubsection.\hskip 2ex}%
 }%
 \begin{#2}}

\renewcommand {\sssbegin}[2][*]
 {\refstepcounter{subsubsection}%\label{sss#1}
\if#1*
\addcontentsline{toc}{subsubsection}{\thesubsubsection.\hskip 1pc #2}%
\else
\addcontentsline{toc}{subsubsection}{\thesubsubsection.\hskip 1pc #2. #1}
\fi
 \def \secno {\gdef \secno {}{\ssecfont \thesubsubsection.\hskip 2ex}%
 }%
 \begin{#2}}

\renewcommand {\parbegin}[2][*]
 {\refstepcounter{paragraph}%\label{ssss#1}
\if#1*
\addcontentsline{toc}{paragraph}{\theparagraph.\hskip 1pc #2}%
\else
\addcontentsline{toc}{paragraph}{\theparagraph.\hskip 1pc #2. #1}
\fi
 \def \secno {\gdef \secno {}{\ssecfont \theparagraph.\hskip 2ex}%
 }%
 \begin{#2}}

\newcommand {\ce}{{\mathrm{CE}}}

\rmnameii{etr}{etr}
\rmnameii{evv}{ev}

\DeclareMathOperator{\K}{\mathbb{K}}
\DeclareMathOperator{\F}{\mathbb{F}}

\DeclareMathOperator{\ab}{\mathrm{ab}}

\newcommand {\w}{\omega}

\newcommand{\fC}{\mathfrak{C}}
\newcommand{\LR}{\mathrm{LR}}
\newcommand{\obs}{\mathrm{obs}}
\newcommand{\fobs}{\mathfrak{obs}}

\begin{document}

\title[Deformations of restricted Lie-Rinehart algebras]{An alternative approach to deformations of restricted Lie-Rinehart algebras}

\author{Quentin Ehret}
\address {Division of Science and Mathematics, New York University Abu Dhabi, P.O. Box 129188, Abu Dhabi, United Arab Emirates.}
\email{qe209@nyu.edu}

\thanks{The author is supported by the grant NYUAD-065.}


\begin{abstract}
In this paper, we introduce a novel approach to the deformation theory of restricted Lie-Rinehart algebras in positive characteristic, based on the deformation theory of restricted morphisms introduced in our earlier work. We provide a full cohomology complex adapted to restricted Lie-Rinehart algebras in characteristic 2, define formal deformations, study obstructions and equivalence classes, and show that these are controlled by the cohomology we introduce. In characteristic $p\geq 3$, we construct 2-cocycles and investigate their relationship with formal deformations. Explicit computations are given in characteristic 2 to illustrate the theory.   

\end{abstract}

%\subjclass[2010]{Primary 17B50; Secondary 17B20}

\maketitle

\thispagestyle{empty}
\setcounter{tocdepth}{2}
\tableofcontents

%%%%%%%%%%%%%%%%%%%%%%%%%%%%%%%%%%%%%%%%%%%%%%%%%%%%%%%%%%%%%%%%%%%%%%%

\section{Introduction} \label{intro}

\textbf{Restricted Lie algebras.} Restricted Lie algebras appeared in the work of Jacobson as a way to axiomatize the properties of derivations of associative algebras that arise only in strictly positive characteristic, see \cite{J}. Roughly speaking, one requires the presence of a so-called \textit{p-map} that behaves like the $p$-th power of associative algebras, see Definition \ref{def:RLie}.  Later, those algebras appeared to be of prime interest due to their links to algebraic groups, and representation theory in positive characteristic (see \cite{SF}). Their cohomology is way more complicated than the usual Chevalley-Eilenberg cohomology for Lie algebras and was first defined by Hochschild in \cite{H}. Explicit formulas in low degrees and interpretation in terms of central extensions were provided by Evans and Fuchs (see \cite{EF}), who also sketched a theory of infinitesimal deformations. The full deformation theory of restricted Lie algebras and its links with the restricted cohomology was done in \cite{EM2}. In this paper, we also exhibit a full cohomology complex in characteristic 2 that have no analogs in other characteristics\footnote{This complex was, in fact, already known to experts.} and show that it fits with formal deformations. Furthermore, we investigate the deformations of restricted morphisms between restricted Lie algebras and introduce a cohomology that controls those deformations, following previous work in characteristic zero (see e.g \cite{Ge, GS1, GS2} for associative algebras, \cite{NR1,NR2} for Lie algebras and \cite{Ma} for Leibniz algebras). Other applications of the restricted cohomology are classification of $p$-nilpotent Lie algebras and superalgebras of small dimension  (see \cite{MS, BE}) and study of extensions of some class of restricted algebras (see e.g. \cite{EFP,EF2,EF3}), to cite a few. Furthermore, the generalization of the restricted cohomology to restricted Lie superalgebras can be found in \cite{BE}. \\

\textbf{Restricted Lie-Rinehart algebras.} Lie-Rinehart algebras over fields of characteristic zero appeared as algebraic analogs of Lie algebroids in differential geometry. The study of these structures was pioneered by Palais (\cite{Pa}), Herz (\cite{He}) and Rinehart (\cite{Ri}). Roughly speaking, a Lie-Rinehart algebra is a triple $(A,L,\rho)$ where $A$ is an associative commutative algebra, $L$ is a Lie algebra and an $A$-module and $\rho:L\rightarrow\Der(A)$ is an $A$-linear morphism of Lie algebras that satisfies a compatibility condition, see Eq. \eqref{eq:leibniz-def}.  Lie-Rinehart algebras have recently gained attention and have been widely studied, especially in the work of Huebschmann (see \cite{Hu1, Hu2, Hu3}). Several generalizations of this concept have since emerged, including Hom-Lie-Rinehart algebras (see \cite{MM1, MM2}), Lie-Rinehart superalgebras (see \cite{C, EM, Ro}), and $n$-Lie-Rinehart (super)algebras (see \cite{BCMS, BCZ}).

In positive characteristic, the Lie algebra is required to be restricted, and additional structures emerge from the interactions between the structural maps of the Lie-Rinehart algebra and the $p$-map of the restricted Lie algebra, see Eq. \eqref{eq:hoch}. Restricted Lie-Rinehart algebras first appeared implicitly in Hochschild's work on his Galois-Jacobson theory of purely inseparable extensions of degree 1 (see \cite{Ho}), as well as in \cite{Ru}. A proper definition was given later by Dokas in \cite{Do}, who notably proved a Jacobson-type theorem for restricted Lie-Rinehart algebras (see \cite[Proposition 2.2]{Do}) and gave a definition of cohomology groups for those objects based on Beck modules and Quillen cohomology, that he used to classify extensions. Moreover, Schauenburg also gave an alternative approach to the notion of restricted universal enveloping algebra of a restricted Lie-Rinehart algebra, see \cite{Sc}.\\

\textbf{Outline of the paper and main results.} This paper is an application of the theory developed in \cite{EM2} to restricted Lie-Rinehart algebras. The aim is to present a new approach to the deformations theory of restricted Lie-Rinehart algebras in positive characteristic, different from the method described in our previous work \cite{EM1}. In this paper, we introduced the notion of \textit{restricted multiderivation} to encode restricted Lie-Rinehart algebras structures. We used this correspondence to build a deformation theory together with a controlling cohomology. We were primarily inspired by methods originating from the study of differential geometric objects, such as Lie algebroids over the real numbers (see \cite{CM, MM2, Ru}). In that case, the anchor map is inherently a Lie morphism, as are all its deformations, provided that the Leibniz rule for Lie algebroids is satisfied. In the purely algebraic case of (restricted) Lie-Rinehart algebras, this property no longer holds. A necessary condition for the anchor map of a Lie-Rinehart algebra $(A, L, \rho)$ to be a Lie morphism is that the Lie algebra $L$ is projective as $A$-module (see e.g. \cite{FMM,Hu1}). However, in the most general case, this is not guaranteed, and additional conditions have to be considered. In \cite{EM2}, we introduced a cohomology theory to control the deformations of restricted morphisms of restricted Lie algebras. In the present paper, we extend those methods from \cite{EM2} to study the deformations of restricted Lie-Rinehart algebras, even when the underlying Lie algebra is not projective as $A$-module. This method appeared to be particularly efficient in characteristic $p=2$.\\

The paper is organized as follows. In Section \ref{sec:recollection}, we review the basics on restricted Lie and Lie-Rinehart algebras. Section \ref{sec:p=2} deals with the characteristic 2 case. We recall the formulas of the restricted cohomology of restricted Lie algebras and restricted morphisms given in \cite{EM2} and we exhibit a well-defined sub-complex adapted to the problem of deformations of restricted Lie-Rinehart algebras, see Theorem \ref{thm:p=2welldefined}. Furthermore, we introduce formal deformations and investigate obstructions to extending finite-order deformations (see Proposition \ref{prop:obs-p=2}), equivalence (see Proposition \ref{prop:equiv-p=2}) and show that they are controlled by the cohomology we have introduced. In Section \ref{sec:ex}, we conclude the study by computing explicitly 2-cocycles on some examples. Next, we move to the characteristic $p\geq3$ case in Section \ref{sec:p-geq3}. The situation here is slightly more complicated and therefore we restrict ourselves to building deformation 2-cocycles that will be used to control formal deformations (see Propositions \ref{prop:cocy-p} and \ref{prop:obs-p}).

%%%%%%%%%%%%%%%%%%%%%%%%%%%%%%%%%%%%%%%%%%%%%%%%%%%%%%%%%%%%%%%%%%%%%%%%%%%%%%%
\section{Recollections on restricted Lie-Rinehart algebras}\label{sec:recollection}
In this section, $\K$ is a field of characteristic $p>0$. We recall the necessary material on restricted Lie algebras and restricted Lie-Rinehart algebras, following \cite{J} and \cite{Do}.
%%%%%%%%%%%%%%%%%%%%%%%%%%%%%%%%%%%
\subsection{Restricted Lie algebras}\label{def:RLie}
Let $L$ be a finite-dimensional Lie algebra over $\K$. Following \cite{J}, a map $(\cdot)^{[p]}:L\rightarrow L, \quad x\mapsto x^{[p]}$ is called a~\textit{$p$-map} (or $p$\textit{-structure}) on $L$ and $L$ is said to be {\it restricted}  if 
\begin{align}\label{def:restrictedLie}
(\lambda x)^{[p]}&=\lambda^p x^{[p]} \text{ for all } x\in L \text{ and for all } \lambda \in \K;\\
\ad_{x^{[p]}}&=(\ad_x)^p \text{ for all } x\in L;\\
(x+y)^{[p]}&=x^{[p]}+y^{[p]}+\displaystyle\sum_{1\leq i\leq p-1}s_i(x,y), \text{ for all } x,y\in L,\label{eq:si}
\end{align}
where the coefficients $s_i(x,y)$ can be obtained from the expansion
\begin{equation}
(\ad_{\lambda x+y})^{p-1}(x)=\displaystyle\sum_{1\leq i \leq p-1} is_i(x,y) \lambda^{i-1}.
\end{equation}
In the case where $p=2$, Condition \eqref{eq:si} reduces to
$$(x+y)^{[2]}=x^{[2]}+y^{[2]}+[x,y],~\forall x,y\in L.$$

\sssbegin{Example}\label{ex:deri_p}
    Let $A$ be an associative algebra. Then, the space of derivations of $A$ is a restricted Lie algebra. Indeed, let
    $$\Der(A):=\{d:A\rightarrow A,~d(ab)=d(a)b+ad(b),~\forall a,b\in A\}.$$ It is well known that $\Der(A)$ is a Lie algebra with the ordinary commutator bracket. Moreover, for all $d\in \Der(A)$ and all $a,b\in A$, we have
    \begin{equation}\label{derip}  d^p(ab)=\sum_{i=0}^{p}\binom{p}{i}d^i(a)d^{p-i}(b)=ad^p(b)+d^p(a)b.   
    \end{equation} Therefore, the map $d\mapsto d^p$ is an endomorphism of $\Der(A)$ and defines a restricted structure on $\Der(A)$.
\end{Example}

\sssbegin{Example}\label{ex:witt}
 Let $\mathbb{K}$ be a field of characteristic $p\geq5$. We consider the \textit{Witt algebra} $W(1)=\text{Span}_{\K} \{e_{-1},e_0,\cdots,e_{p-2}\}$   with the bracket
		%\begin{center}	
		$$[e_i,e_j]:=
		\begin{cases}
			(j-i)e_{i+j} ~\text{ if }~ i+j\in\{-1,...,p-2\};\\
			0 ~\text{ otherwise,}\\
		\end{cases}$$
		%\end{center}
		and the $p$-map
		%\begin{center}	
		$$e_i^{[p]}:=
		\begin{cases}
			e_0~\text{ if } i=0;\\
			0 ~~\text{  if } i\neq 0.\\
		\end{cases}$$
		%\end{center}
		
		Then $\left(W(1),[\cdot,\cdot], (\cdot)^{[p]}\right)$ is a restricted Lie algebra (see \cite{EFP}). More examples can be found in \cite{Ko} and in \cite{MS}.
\end{Example}

\subsubsection{Morphisms and representations} Let $\left( L,[\cdot,\cdot]_L,(\cdot)^{[p]_L}\right) $ and $\left( H,[\cdot,\cdot]_H,(\cdot)^{[p]_H}\right) $ be two restricted Lie algebras. A Lie algebra morphism $\varphi:L\longrightarrow H$ is called \textit{restricted morphism} (or \textit{$p$-morphism}) if it satisfies \begin{equation}
\varphi\bigl(x^{[p]_L} \bigl)=\varphi(x)^{[p]_H}\quad\forall x\in L.
\end{equation}
An $L$-module $M$ is called \textit{restricted} if in addition to the ordinary condition
 \begin{equation}
   [x,y]_L\cdot m=x\cdot(y\cdot m)-y\cdot(x\cdot m) ,~\forall x,y\in L,\forall m\in M,
\end{equation} it also satisfies
\begin{equation}
\underbrace{x\cdots x}_{p\text{~~terms}}\cdot m =x^{[p]_L}\cdot m, \quad  \forall x\in L,~\forall  m\in M.
\end{equation}
A linear map $d:L\rightarrow L$ is called a \textit{restricted derivation} of $L$  if in addition to the ordinary condition
 \begin{equation}
    d\bigl([x,y]_L\bigl)=[d(x),y]_L+[x,d(y)]_L,~\forall x,y\in L,
\end{equation} it also satisfies
\begin{equation}
    d\bigl(x^{[p]_L}\bigl)=\ad_x^{p-1}\circ d(x),~\forall x\in L.
\end{equation}

The following theorem, due to Jacobson, is useful to investigate $p$-structures on a Lie algebra.
\sssbegin{Theorem}[\cite{J}] \label{Jac}
Let $(e_j)_{j\in J}$ be a~basis of $L$ such that there are $f_j\in L$ satisfying $(\ad_{e_j})^p=\ad_{f_j}$. Then, there exists exactly one $p$-map  $(\cdot)^{[p]}:L\rightarrow L$ such that 
\[
e_j^{[p]}=f_j \quad \text{ for all $j\in J$}.
\]
\end{Theorem}

%%%%%%%%%%%%%%%%%%%%%%%%%%%%%%%%%%%%%%%%%%%%%%%%%%%%%%%%%%%%%%%%%%%%%%%%

\subsection{Lie-Rinehart algebras}
%%%%%%%%%%%%%%%%%%%%%%%%%%%%%%%%%%%%%%%%%%%%%%%%%%%%%%%%%%%%%%%%%%%%%%%%%%%%%%
Following \cite{Hu1,Ri}, a \textit{Lie-Rinehart algebra} is given by a triple $(A,L,\rho)$, where $A$ is a commutative associative algebra, $(L,[\cdot,\cdot])$ is a Lie algebra and an $A$-module, and $\rho:L\rightarrow\Der(A)$ is an $A$-linear morphism of Lie algebras satisfying
\begin{equation}\label{eq:leibniz-def}
    [x,ay]=a[x,y]+\rho(x)(a)y,~\forall x,y\in L,~\forall a\in A.
\end{equation}

\sssbegin{Example}\label{ex:der} Let $A$ be a commutative associative algebra. Then, $(A,\Der(A),\id)$ is a Lie-Rinehart algebra.
\end{Example}

A \textit{representation} of a Lie-Rinehart algebra $(A,L,\rho)$ is an $A$-module $M$ endowed with a $A$-linear Lie algebra morphism $\pi:L\rightarrow \text{End}(M)$ satisfying \begin{equation} \pi(x)(am)=a\pi(x)(m)+\rho(x)(a)m,~~~ \forall a\in A,~\forall m\in M,~ \forall x\in L.  \end{equation}

In positive characteristic, the following Lemma holds, which is a reformulation of \cite[Lemma 1]{Ho}.

\sssbegin{Lemma}[\text{see \cite[Lemma 2.1]{Sc}}]\label{lem:sc}
	Let $M$ be a representation of a Lie-Rinehart algebra $(A,L,\rho)$. Then, we have 
    \begin{equation}    \pi(ax)^p=a^p\pi(x)^p+\rho(ax)^{p-1}(a)\pi_x,~\forall x\in L,~\forall a\in A.   
    \end{equation}
\end{Lemma}
Applied to the Lie-Rinehart algebra $(A,\Der(A),\id)$ (see Example \ref{ex:der})  this Lemma gives 
 \begin{equation} (aD)^p=a^pD^p+(aD)^{p-1}(a)D,~\forall a\in A,~\forall D\in\Der(A).   \end{equation}
 
 This lead to the following general definition of a restricted Lie-Rinehart algebra.

 \subsection{Restricted Lie-Rinehart algebras}\label{def:RLRp} (see \cite{Do}) A \textit{restricted Lie-Rinehart algebra} is given by a triple $(A,L,\rho)$, where $A$ is a commutative associative algebra, $\bigl(L,[\cdot,\cdot],(\cdot)^{[p]}\bigl)$ is a restricted Lie algebra, and $\rho:L\rightarrow\Der(A)$ is an $A$-linear restricted Lie morphism satisfying
    \begin{align}\label{eq:leib}
        [x,ay]&=a[x,y]+\rho(x)(a)y,~\forall x,y\in L,~\forall a\in A; \end{align}
        and
        \begin{align}\label{eq:hoch}
        (ax)^{[p]}&=a^px^{[p]}+\rho(ax)^{p-1}(a)x, ~\forall x\in L,~\forall a\in A.
    \end{align}
    
\sssbegin{Example}
    Let $A$ be an associative commutative algebra and let $(A,\Der(A),\id)$ be the Lie-Rinehart algebra of Example \ref{ex:der}. Then, this Lie-Rinehart algebra is restricted. This applies in particular to the restricted Witt algebra (see Example \ref{ex:witt}), which is the derivations algebra of the the commutative associative algebra  $A:=\K[x]/(x^p-1)$, see \cite{EF02} for more details.
\end{Example}

\sssbegin{Example}
    Let $(A,\{\cdot,\cdot\},)$ be a Poisson algebra with unit $1_A$ equipped with a $p$-map $(\cdot)^{\{p\}}$ (that is, $(A,\{\cdot,\cdot\},(\cdot)^{\{p\}})$ is a restricted Lie algebra). Following \cite{PS,BYZ}, $(A,\{\cdot,\cdot\},(\cdot)^{\{p\}})$ is called \textit{restricted Poisson algebra} if
    \begin{align}
        (xy)^{\{2\}}&=x^2y^{\{2\}}+ y^2x^{\{2\}}+xy\{x,y\},&\forall x,y\in A,\quad \text{if }p=2;\\
        (x^2)^{\{p\}}&=2x^px^{\{p\}},&\forall x,y\in A,\quad \text{if }p\geq3.
    \end{align}
 Consider the module of Kähler differentials $\Omega^1(A)$ of a restricted Poisson algebra $(A,\{\cdot,\cdot\},(\cdot)^{\{p\}})$, which is the free $A$-module generated by symbols $du,~u\in A,$ with relations
    \begin{equation}
        d(u+v)=du+dv,~d(uv)=udv+vdu,~d1_A=0.
    \end{equation}
    Then, it has been shown in \cite[Theorem 3.8]{Hu1} that $(A,\Omega^1(A),\rho)$ is a Lie-Rinehart algebra, with
    \begin{align}
     [xdu,ydv]_{\Omega^1(A)}&:=x\{u,y\}d v+y\{x,v\}d u+xyd\{u,v\},\quad &\forall x,y,u,v\in A;\\
     \rho(xdu)&:=x\{u,-\} &\forall x,u\in A.
    \end{align}
    Suppose that $\Omega^1(A)$ is free as an $A$-module. Then, $(A,\Omega^1(A),\rho)$ becomes a restricted Lie-Rinehart algebra, with the $p$-map (see \cite[Theorem 4.3.2]{BEL} and \cite[Theorem 8.2]{BYZ}) 
    \begin{equation}
        (xdu)^{[p]}:=x^pd(u^{\{p\}})+\rho(xdu)^{p-1}(x)du,\quad \forall x,u\in A.
    \end{equation}
\end{Example}

Let $(A,L,\rho)$ and $(A,\tilde{L},\tilde{\rho})$ be restricted Lie-Rinehart algebras sharing the same associative commutative algebra $A$. A morphism of restricted Lie-Rinehart algebras is an $A$-linear restricted Lie morphism $\phi:L\rightarrow\tilde{L}$ satisfying
\begin{equation}\label{eq:morphRLR}
    \rho(x)=\tilde{\rho}\circ\phi(x),~\forall x\in L.
\end{equation}

\subsection{Derivations of the formal space $A[[t]]$} Let $A$ be an associative commutative algebra. The formal space $A[[t]]$ is defined by
\begin{equation}
    A[[t]]:=\Bigl\{\overline{a}:=\sum_{i\geq 0}t^ia_i,~a_i \in A\Bigl\}.
\end{equation}
The multiplication on $A[[t]]$ is given by the Cauchy product, that is,
\begin{equation}
    \overline{a}\overline{b}=\Bigl(\sum_{i\geq 0}t^ia_i\Bigl)\Bigl(\sum_{j\geq 0}t^jb_j\Bigl):=\sum_{i\geq 0}t^i\sum_{k=0}^{i}a_k b_{i-k},~~\forall~\overline{a},\overline{b}\in A[[t]].
\end{equation}

\sssbegin{Lemma}
Let $D\in\Der(A)$. Then, $D$ extends to an $A[[t]]$-linear derivation of $A[[t]]$.
\end{Lemma}
\begin{proof}
Let $\overline{a},\overline{b}\in A[[t]]$. We have
\begin{align*}
    D\bigl(\overline{a}\overline{b}\bigl)
    &=D\Bigl(\sum_{i\geq 0}t^i\sum_{k= 0}^{i}a_k b_{i-k}\Bigl)=\sum_{i\geq 0}t^i\sum_{k= 0}^{i}\bigl(D(a_k)b_{i-k}+a_kD(b_{i-k})\bigl)\\
    &=\Bigl(\sum_{i\geq0}t^iD(a_i)\Bigl)\Bigl(\sum_{j\geq0}t^jb_j\Bigl)+\Bigl(\sum_{i\geq0}t^ia_i\Bigl)\Bigl(\sum_{j\geq0}t^jD(b_j)\Bigl)\\
    &=D\Bigl(\sum_{i\geq0}t^ia_i\Bigl)\Bigl(\sum_{j\geq0}t^jb_j\Bigl)+\Bigl(\sum_{i\geq0}t^ia_i\Bigl)D\Bigl(\sum_{j\geq0}t^jb_j\Bigl)\\
    &=D\bigl(\overline{a}\bigl)\overline{b}+\overline{a}D\bigl(\overline{b}\bigl).
\end{align*}

\end{proof}

%%%%%%%%%%%%%%%%%%%%%%%%%%%%%%%%%%%%%%%%%%%%%%%%%%%%%%%%%%%%%%%%%%%%%%%%%%%%%%%
\section{Deformations of restricted Lie-Rinehart algebras, $p=2$}\label{sec:p=2} In this section, $\K$ is a field of characteristic $p=2$. In this particular case, the Definition \ref{def:RLRp} of a restricted Lie-Rinehart algebra $(A,L,\rho)$ becomes much simpler. Indeed, Eq. \eqref{eq:hoch} can be stated as
    \begin{equation}\label{eq:hoch2} 
            (ax)^{[2]}=a^2x^{[2]}+\rho(ax)(a)x, ~\forall x\in L,~\forall a\in A.
    \end{equation}
    This specificity allows to build a full deformation cohomology complex, inspired by the deformation cohomology of restricted morphisms introduced in \cite{EM2}. We construct this complex, show that it is well-defined (see Theorem \ref{thm:p=2welldefined}), develop the theory of formal deformations of restricted Lie-Rinehart algebras in characteristic 2 (see Section \ref{sec:deformations}) and compute some examples (see Section \ref{sec:ex}).

 \subsection{Restricted cohomology of restricted Lie algebras, $p=2$} In this section, we recall the construction of the restricted cohomology for restricted Lie algebras over a field of characteristic $p=2$, that we described in \cite{EM2}. 
	
\subsubsection{Restricted cochains, $p=2$}	
Let $\bigl(L,[\cdot,\cdot],(\cdot)^{[2]}\bigl)$ be a restricted Lie algebra and let $M$ be a restricted $L$-module. We set $C_{\mathrm{res}}^0(L,M):=C_{\mathrm{CE}}^0 (L,M)$ and $C_{\mathrm{res}}^1(L,M):=C_{\mathrm{CE}}^1 (L,M)$\footnote{The subscripts ``CE" refer to the ordinary Chevalley-Eilenberg cohomology, see \cite{EM2}.}.

Let $n\geq2$, $\varphi\in C_{\mathrm{CE}}^n (L,M)$, $\omega:L\times \wedge^{n-2}L\rightarrow M$, $\lambda\in\F$ and $x,z_2,\cdots,z_{n-1}\in L$. The pair $(\varphi,\omega)$ is a $n$-cochain of the restricted cohomology if
		\begin{align}
			\omega(\lambda x, z_2,\cdots,z_{n-1})&=\lambda^2\omega(x,z_2,\cdots,z_{n-1});\\
			\omega(x+y,z_2,\cdots,z_{n-1})&=\omega(x,z_2,\cdots,z_{n-1})+\omega(y,z_2,\cdots,z_{n-1})\\\nonumber&~\,\,+\varphi(x,y,z_2,\cdots,z_{n-1});\\
            \text{the map } (z_2,\cdots,z_{n-1})&\mapsto\omega(\cdot,z_2,\cdots,z_{n-1}) \text{ is multilinear.}
		\end{align}
		We denote by $C_{\mathrm{res}}^n(L;M)$ the space of $n$-cochains of $L$ with values in $M$.

\subsubsection{Restricted differential operators, $p=2$}
For $n\geq 2$, the differential maps\\
	$d^n_{\mathrm{res}}:C^n_{\mathrm{res}}(L;M)\rightarrow C^{n+1}_{\mathrm{res}}(L;M)$ are given by $d^n_{\mathrm{res}}(\varphi,\omega)=\bigl(d^n_{\mathrm{CE}}(\varphi),\delta^n(\omega)\bigl),$ where
	\begin{align*}	\delta^n\omega(x,z_2,\cdots,z_n)&:=x\cdot\varphi(x,z_2,\cdots,z_n)+\sum_{i=2}^{n}z_i\cdot\omega(x,z_2,\cdots,\hat{z_i},\cdots,z_n)\\	&+\varphi(x^{[2]},z_2,\cdots,z_n)+\sum_{i=2}^{n}\varphi\left([x,z_i],x,z_2,\cdots,\hat{z_i},\cdots,z_n \right)\\
		&+\sum_{2\leq i<j\leq n}\omega\left(x,[z_i,z_j],z_2,\cdots,\hat{z_i},\cdots,\hat{z_j},\cdots,z_n  \right).
	\end{align*}
    For $n=0,1$, we define $d^0_{\mathrm{res}}=d^0_{\mathrm{CE}}$ and 
	\begin{align*}
		d^1_{\mathrm{res}}:C_{\mathrm{res}}^1(L;M)&\rightarrow C_{\mathrm{res}}^2(L;M)\\
		\varphi&\mapsto\bigl(d_{\mathrm{CE}}^1\varphi,\delta^1\varphi\bigl),~\text{where }\delta^1\varphi(x):=\varphi\bigl( x^{[2]} \bigl)+x\cdot\varphi(x),~\forall x\in L.
	\end{align*}

	\sssbegin{Theorem}[\text{\cite[Theorem 4.9]{EM2}}]\label{thm:cohomology2}
		Let $\bigl(L,[\cdot,\cdot],(\cdot)^{[2]}\bigl)$ be a restricted Lie algebra and let $M$ a restricted $L$-module. The complex $\left(C_{\mathrm{res}}^n(L;M),d^n_{\mathrm{res}} \right)_{n\geq 0}$ is a cochain complex. The $n^{th}$ restricted cohomology group of the Lie algebra $L$ in characteristic 2 is defined by	
		$$ H_{\mathrm{res}}^n(L;M):=Z_{\mathrm{res}}^n(L;M)/B^n_{\mathrm{res}}(L;M),$$ 
		with $Z_{\mathrm{res}}^n(L;M)=\Ker(d^n_{\mathrm{res}})$ the restricted $n$-cocycles and $B_{\mathrm{res}}^n(L;M)=\text{Im}(d^{n-1}_{\mathrm{res}})$ the restricted $n$-coboundaries.
	\end{Theorem}

\subsubsection{Restricted cohomology of restricted morphisms, $p=2$}   
Let $(L,[ \cdot , \cdot ]_L,(\cdot )^{[2]_L})$ and $(M,[ \cdot , \cdot ]_M,(\cdot )^{[2]_M})$ be restricted Lie algebras, let $\varphi:L\rightarrow M$ be a restricted morphism and let $n\geq1$. We define 
	$$\fC^n_{\mathrm{res}}(\varphi,\varphi):=C^n_{\mathrm{res}}(L;L)\times  C^n_{\mathrm{res}}(M;M)\times C^{n-1}_{\mathrm{res}}(L;M),$$
	and $\fC^0_{\mathrm{res}}(\varphi,\varphi):=0$. For $n\geq 3$, the differential maps are given by
	\begin{align}
		\fd^n_{\mathrm{res}}:\fC^n_{\mathrm{res}}(\varphi,\varphi)&\rightarrow\fC^{n+1}_{\mathrm{res}}(\varphi,\varphi)\nonumber\\\label{diffmorph2}
    \Bigl((\mu,\w),(\nu,\epsilon),(\theta,\gamma)\Bigl)&\mapsto\Bigl((d^n_{\ce}\mu,\delta^n\w),(d^n_{\ce}\nu,\delta^n\epsilon),\bigl(\alpha_{\mu,\nu}(\theta),\beta_{\w,\epsilon}(\gamma)\bigl) \Bigl),    
		% \begin{pmatrix}(\mu,\w)\\(\nu,\epsilon)\\(\theta,\gamma)\end{pmatrix}&\mapsto \begin{pmatrix}\bigl(d^n_{\ce}\mu,\delta^n\w\bigl)\\\bigl(d^n_{\ce}\nu,\delta^n\epsilon\bigl)\\\bigl(\alpha_{\mu,\nu}(\theta),\beta_{\w,\epsilon}(\gamma)\bigl)\end{pmatrix},
	\end{align}
	where \begin{equation}\alpha_{\mu,\nu}(\theta):=\varphi\circ\mu+\nu\circ\varphi^{\otimes n}+d^{n-1}_{\ce}\theta\end{equation} and \begin{equation}\beta_{\w,\epsilon}(\gamma):=\varphi\circ\w+\epsilon\circ \varphi^{\otimes(n-1)}+\delta^{n-1}\gamma.\end{equation} 
    Moreover, we have
    \begin{align}
        \fd_{\rm res}^1:\fC^1_{\mathrm{res}}(\varphi,\varphi)&\rightarrow\fC^{2}_{\mathrm{res}}(\varphi,\varphi)\\\nonumber
        (\mu,\nu,m)&\mapsto \Bigl((d^1_{\ce}\mu,\delta^1\mu),(d^1_{\ce}\nu,\delta^1\nu),\alpha_{\mu,\nu}(m)  \Bigl)
    \end{align} and
    \begin{align}
        \fd_{\rm res}^2:\fC^2_{\mathrm{res}}(\varphi,\varphi)&\rightarrow\fC^{3}_{\mathrm{res}}(\varphi,\varphi)\\\nonumber
        \Bigl((\mu,\w),(\nu,\epsilon),\theta\Bigl)&\mapsto \Bigl((d^2_{\ce}\mu,\delta^2\w),(d^2_{\ce}\nu,\delta^2\epsilon),\bigl(\alpha_{\mu,\nu}(\theta),\beta_{\mu,\nu}(\theta)\bigl)  \Bigl).
    \end{align}    
    
     We denote by $\mathfrak{Z}^n_{\mathrm{res}}(\varphi,\varphi):=\Ker(\fd^n_{\mathrm{res}})$ and $\mathfrak{B}^n_{\mathrm{res}}(\varphi,\varphi):=\text{Im}(\fd^{n-1}_{\mathrm{res}})$.
    \sssbegin{Theorem}[\text{\cite[Theorem 4.22]{EM2}}]\label{thm:cohomorph2}
    For all $n\in\mathbb{N}$, the maps $\fd^n_{\mathrm{res}}$ are well defined and satisfy $\fd^{n+1}_{\mathrm{res}}\circ\fd^n_{\mathrm{res}}=0$.
\end{Theorem}

\subsection{Deformation cohomology of restricted Lie-Rinehart algebras, $p=2$} In this section, we introduce a cohomology of restricted Lie algebras derived from the restricted morphisms cohomology, which is different from the cohomology given in \cite{BEL}. Let $(A,L,\rho)$ be a restricted Lie-Rinehart algebra and let $M$ be a restricted Lie-Rinehart module. For every derivation $d\in \Der(A)$, we define 
\begin{equation}
     C^1_{\LR}(L;M)^d:=\bigl\{(\mu,d),~\mu\in C_{\mathrm{res}}^1(L;M),~\mu(ax)=a\mu(x)+d(a)x,\forall x\in L,\forall a\in A  \bigl\}.
\end{equation}
For every $A$-linear map $\theta: L\rightarrow\Der(A)$, we define 
\begin{equation}
    C^2_{\LR}(L;M)^{\theta}:=\bigl\{(\mu,\w,\theta),~(\mu,\w)\in C_{\mathrm{res}}^2(L;M), \text{ satisfying } \eqref{eq:LRcohomop=2a} \text{ and } \eqref{eq:LRcohomop=2b}  \bigl\},
    \end{equation}
where
\begin{align}\label{eq:LRcohomop=2a}
\mu(x,ay)&=a\mu(x,y)+\theta(x)(a)y,~\forall a\in A,\forall x,y\in L~\text{ and}\\\label{eq:LRcohomop=2b}
\w(ax)&=a^2\w(x)+\theta(ax)(a)x,~\forall a\in A,\forall x\in L.
\end{align}

Let $n\geq 3$. For every $A$-linear maps $(\theta,\gamma)\in C_{\mathrm{res}}^{n-1}\bigl(L;\Der(A)\bigl)$, we define
\begin{equation}
    C^n_{\LR}(L;M)^{(\theta,\gamma)}:=\bigl\{(\mu,\w,\theta,\gamma),~(\mu,\w)\in C_{\mathrm{res}}^n(L;M), \text{ satisfying } \eqref{eq:LRcohomop=2c},~\eqref{eq:LRcohomop=2d} \text{ and } \eqref{eq:LRcohomop=2e}  \bigl\},
    \end{equation}
where, for all $x,z_1,\cdots,z_{n-2}\in L$ and all $a\in A$,
\begin{align}\label{eq:LRcohomop=2c}
    \mu(x_1,\cdots,x_{n-1},ax_n)&=a\mu(x_1,\cdots,x_n)+\theta(x_1,\cdots,x_{n-1})(a)x_n;\\\label{eq:LRcohomop=2d}
    \w(ax,z_1,\cdots,z_{n-2})&=a^2\w(x,z_1,\cdots,z_{n-2})+\theta(ax,z_1,\cdots,z_{n-2})(a)x;\\\label{eq:LRcohomop=2e}
    \w(x,z_1,\cdots,az_i,\cdots,z_{n-2})&=a\w(z_1,\cdots,z_i,\cdots,z_{z-2})+\gamma(x,z_1,\cdots,\hat{z_i},\cdots,z_n)(a)z_i.\end{align}

We denote by
\begin{equation}
    C^1_{\LR}(L;L):=\bigcup_{d}C^1_{\LR}(L;L)^{d};~C^2_{\LR}(L;L):=\bigcup_{\theta}C^2_{\LR}(L;L)^{\theta};
\end{equation}
and
\begin{equation}
    C^n_{\LR}(L;L):=\bigcup_{(\theta,\gamma)}C^n_{\LR}(L;L)^{(\theta,\gamma)},~\text{for all } n\geq 3.
\end{equation}
In the case where $M=L$, the space $C^n_{\LR}(L;L),~n\geq 1$, is called the $n^{\text{th}}$ \textit{restricted Lie-Rinehart deformation cohomology cochains space}. There are natural embeddings
\begin{align}
    \iota_1:& ~C^1_{\LR}(L;L)^{d}~\hookrightarrow\fC^1_{\mathrm{res}}(\rho,\rho),\quad (\mu,d)\mapsto (\mu,0,d );\\
    \iota_2:& ~C^2_{\LR}(L;L)^{\theta}~\hookrightarrow\fC^2_{\mathrm{res}}(\rho,\rho),\quad (\mu,\w,\theta)\mapsto \Bigl((\mu,\w),(0,0),\theta \Bigl);\\
    \iota_n:& ~C^n_{\LR}(L;L)^{(\theta,\gamma)}\hookrightarrow\fC^n_{\mathrm{res}}(\rho,\rho),\quad (\mu,\w,\theta,\gamma)\mapsto \Bigl((\mu,\w),(0,0),(\theta,\gamma) \Bigl),~\forall n\geq 3.
\end{align}

The next three Lemmas aim to build differential maps, using the maps given in Eq. \eqref{diffmorph2}. Since the spaces and maps of order 1 and 2 have a slightly different behavior compared to higher order spaces and maps, we consider them separately. In the sequel, whenever a symbol $\iota^{-1}_k(X),~k\geq 1,~X\subset \fC^k_{\mathrm{res}}(\rho,\rho)$ appears, it should be understood as the preimage of $X$ under $\iota_k$.

\sssbegin{Lemma}\label{lem:diff1}
    Let $(A,L,\rho)$ be a restricted Lie-Rinehart algebra and let $d\in\Der(A)$. Then, we have \begin{equation}\iota_2^{-1}\circ\fd^1_{\mathrm{res}}\circ\iota_1\Bigl(C^1_{\LR}(L;L)^d\Bigl)\subset C^2_{\LR}(L;L).\end{equation} 
\end{Lemma}

\begin{proof}
    Let $(\mu,d)\in C^1_{\LR}(L;L)^{d}$. We have
    \begin{equation}    \iota_2^{-1}\circ\fd^1_{\mathrm{res}}\circ\iota_1(\mu,d)=\Bigl(d^1_{\ce}\mu,\delta^1\mu,\alpha_{\mu,0}(d)\Bigl).
    \end{equation} Let $a\in A$ and $x,y\in L$. We have
    \begin{align*}
        d^1_{\ce}\mu(x,ay)=&\mu\bigl([x,ay]\bigl)+[\mu(x),ay]+[x,\mu(ay)]\\
        =&\mu\bigl(a[x,y]\bigl)+\mu\bigl(\rho(x)(a)y\bigl)+a[\mu(x),y]+\rho(\mu(x))(a)y+[x,a\mu(y)]+[x,d(a)y]\\
        =&\underline{a\mu\bigl([x,y]\bigl)}+d(a)[x,y]+\rho(x)(a)\mu(y)+\overline{d\circ\rho(x)(a)y}\\
        &+\underline{a[\mu(x),y]}+\overline{\rho(\mu(x))(a)y}+\underline{a[x,\mu(y)]}+\rho(x)(a)\mu(y)\\
        &+d(a)[x,y]+\overline{\rho(x)\circ d(a)y}\\
        =&ad^1_{\ce}\mu(x,y)+\alpha_{\mu,0}(d)(a)y,
    \end{align*} where the underlined terms give $ad^1_{\ce}\mu(x,y)$ and the overlined terms give $\alpha_{\mu,0}(d)(x)(a)y$. Thus, the pair $\bigl(d^1_{\ce}\mu,\alpha_{\mu,0}(d)\bigl)$ satisfy Eq. \eqref{eq:LRcohomop=2a}. Moreover, $\alpha_{\mu,0}(d)(x)\in\Der(A)$ for all $x\in L$ and we have
    \begin{align*}\nonumber
        \alpha_{\mu,0}(d)(ax)(a)x&=\rho\circ\mu(ax)(a)x+\rho(ax)\circ d(a)x+d\circ\rho(ax)(a)x\\\nonumber
        &=\rho\bigl(a\mu(x)\bigl)(a)x+\rho\bigl(d(a)x\bigl)(a)x+\rho(ax)\circ d(a)x+d\circ\rho(ax)(a)x\\\nonumber
        &=a\rho\bigl(\mu(x)\bigl)(a)x+d(a)\rho(x)(a)x+a\rho(x)\bigl(d(a)\bigl)x+d\bigl(\rho(ax)(a)\bigl)x\\
        &=a\rho\bigl(\mu(x)\bigl)(a)x+a\rho(x)\circ d(a)x+ad\circ\rho(x)(a)x.
    \end{align*}
    Furthermore, we have
    \begin{align*}
        \delta^1(ax)&=[ax,\mu(ax)]+\mu\bigl((ax)^{[2]}\bigl)\\
        &=a[x,\mu(ax)]+\rho\bigl(\mu(ax)\bigl)(a)x+\mu\bigl(a^2x^{[2]}\bigl)+\mu\bigl(\rho(ax)(a)x\bigl)\\
        &=a[x,a\mu(x)]+a[x,d(a)x]+\rho\bigl(a\mu(x)\bigl)(a)x+\rho\bigl(d(a)x\bigl)\\
        &~+a^2\mu(x^{[2]})+\underset{=~0}{\underbrace{d(a^2)x^{[2]}}}+\rho(ax)(a)\mu(x)+d\circ\rho(ax)(a)x\\
        &=a^2[x,\mu(y)]+a\rho(x)(a)\mu(x)+a\rho(x)\circ d(a)x+a\rho\bigl(\mu(x)\bigl)(a)x+d(a)\rho(x)(a)x\\
        &~+a^2\mu(x^{[2]}+a\rho(x)(a)\mu(x)+d(a)\rho(x)(a)x+ad\circ\rho(x)(a)x\\
        &=a^2\delta^1\mu(x)+a\rho(x)\circ d(a)x+a\rho\bigl(\mu(x)\bigl)(a)x+d\circ\rho(x)(a)x\\
        &=a^2\delta^1\mu(x)+\alpha_{\mu,0}(d)(ax)(a)x.
    \end{align*} Therefore, Eq. \eqref{eq:LRcohomop=2b} is satisfied and $\iota_2^{-1}\circ\fd^1_{\mathrm{res}}\circ\iota_1(\mu,d)\in  C^2_{\LR}(L;L)^{\alpha_{\mu,0}(d)}$.
\end{proof}

\sssbegin{Lemma}\label{lem:diff2}
Let $(A,L,\rho)$ be a restricted Lie-Rinehart algebra and let $\theta: L\rightarrow\Der(A)$ an $A$-linear map. Then, we have
\begin{equation}\iota_3^{-1}\circ\fd^2_{\mathrm{res}}\circ\iota_2\Bigl(C^2_{\LR}(L;L)^{\theta}\Bigl)\subset C^3_{\LR}(L;L).\end{equation}
\end{Lemma}

\begin{proof}
Let $(\mu,\w,\theta)$. We have \begin{equation}    \iota_3^{-1}\circ\fd^2_{\mathrm{res}}\circ\iota_2(\mu,\w,\theta)=\Bigl(d^2_{\ce}\mu,\delta^2\mu,\alpha_{\mu,0}(\theta),\beta_{\w,0}(\theta)\Bigl).
    \end{equation}
    We show that $\bigl(\delta^2\w,\beta_{\w,0}(\theta)\bigl)$ satisfy Eq. \eqref{eq:LRcohomop=2e}. Let $x,z\in L$ and $a\in A$. We have
   % \vspace{-2cm}
    \begin{align*}
        \delta^2\w(x,az)&=\mu(x^{[2]},az)+[\w(x),az]+[\mu(x,az),x]+\mu([x,z],x)\\
        &=a\mu(x^{[2]},z)+\theta(x^{[2]})(a)z+a[\w(x),z]+\rho\circ\w(x)(a)z\\
        &~~+\bigl[a\mu(x,z)+\theta(x)(a)z,x\bigl]+\mu\bigl(a[x,z]+\rho(x)(a)z,x\bigl)\\
        &=\underline{a\mu(x^{[2]},z)}+\theta(x^{[2]})(a)z+\underline{a[\w(x),z]}+\rho\circ\w(x)(a)z\\
        &~~+\underline{a[\mu(x,z),x]}+\rho(x)(a)\mu(x,z)+\theta(x)(a)[z,x]+\rho(x)\circ\theta(x)(a)z\\
        &~~+\underline{a\mu([x,z],x)}+\theta(x)(a)[x,z]+\rho(x)(a)\mu(z,x)+\theta(x)\circ\rho(x)(a)z\\
        &=a\delta^2\w(x,z)+\beta_{\w,0}(\theta)(x)(a)z,
    \end{align*} where the underlined terms corresponds to $a\delta^2\w(x,z)$. Similarly, we show that $\bigl(d^2_{\ce}\mu,\alpha_{\mu,0}(\theta)\bigl)$ satisfy Eq. \eqref{eq:LRcohomop=2c} and $\bigl(\delta^2\w,\alpha_{\mu,0}(\theta)\bigl)$ satisfy Eq. \eqref{eq:LRcohomop=2d}. Moreover, $\alpha_{\mu,0}(\theta)(x,y)\in \Der(A)$ and $\beta_{\w,0}(\theta)(x)\in \Der(A)$ for all $x,y\in L$.
\end{proof}

\sssbegin{Lemma}\label{lem:diffn}
For $n\geq3$, let $(A,L,\rho)$ be a restricted Lie-Rinehart algebra and let $(\theta,\gamma)\in C_{\mathrm{res}}^{n-1}(L;L)$. Then, we have
\begin{equation}\iota_{n+1}^{-1}\circ\fd^n_{\mathrm{res}}\circ\iota_n\Bigl(C^n_{\LR}(L;L)^{(\theta,\gamma)}\Bigl)\subset C^{n+1}_{\LR}(L;L).\end{equation}
\end{Lemma}

\begin{proof}
    Let $(\mu,\w,\theta,\gamma)\in C^n_{\LR}(L;L)^{(\theta,\gamma)}$. We have
\begin{equation}\iota_{n+1}^{-1}\circ\fd^n_{\mathrm{res}}\circ\iota_n(\mu,\w,\theta,\gamma)=\bigl(d^n_{\ce}\mu,\delta^n\w,\alpha_{\mu,0}(\theta),\beta_{\w,0}(\gamma)\bigl).\end{equation}
Let $x_1,\cdots,x_{n+1}\in L$ and $a\in A$. We denote by $X:=(x_1,\cdots,x_n)$, $\widehat{X_i}:=(x_1,\cdots,\hat{x_i},\cdots,x_n)$ and $\widehat{X_{ij}}:=(x_1,\cdots,\hat{x_i},\cdots,\hat{x_j},\cdots,x_n)$, where the hat means that the term is omitted. We have
\begin{align*}
    &d_{\ce}^n\mu(X,ax_{n+1})\\&=\sum_{1\leq i<j\leq n}\mu\bigl([x_i,x_j],x_1,\widehat{X_{ij}},ax_{n+1}\bigl)+\sum_{i=1}^{n}\mu\bigl([x_i,ax_{n+1}],\widehat{X_i}\bigl)+\sum_{i=1}^{n}\bigl[x_i,\mu(\widehat{X_i},ax_{n+1})\bigl] +\bigl[ax_{n+1},\mu(X)\bigl]\\
&=\sum_{1\leq i<j\leq n}a\mu\bigl([x_i,x_j],\widehat{X_{ij}},x_{n+1} \bigl)+
\sum_{i=1}^n\theta\bigl([x_i,x_j],\widehat{X_{ij}}\bigl)(a)x_{n+1}+\sum_{i=1}^n\mu\bigl(a[x_i,x_{n+1}],\widehat{X_i}\bigl)\\
&~+\sum_{i=1}^n\mu\bigl(\rho(x_i)(a)x_{n+1},\widehat{X_i}\bigl)+\sum_{i=1}^n\bigl[x_i,a\mu(\widehat{X_i},x_{n+1})  \bigl]+\sum_{i=1}^n\bigl[x_i,\theta(\widehat{X_i})(a)x_{n+1}   \bigl]\\
&~+\rho\circ\mu(X)(a)x_{n+1}+a\bigl[x_{n+1},\mu(X)\bigl]\\
&=\underline{a\sum_{1\leq i<j\leq n}\mu\bigl([x_i,x_j],\widehat{X_{ij}},x_{n+1} \bigl)}+\sum_{i=1}^n\theta\bigl([x_i,x_j],\widehat{X_{ij}}\bigl)(a)x_{n+1}+\underline{a\sum_{i=1}^n\mu\bigl([x_i,x_{n+1}],\widehat{X_i}\bigl)}\\
&~+\sum_{i=1}^n\theta(\widehat{X_i})(a)[x_i,x_{n+1}]+\rho(x_i)(a)\sum_{i=1}^n\mu\bigl(x_{n+1},\widehat{X_i}\bigl)+\sum_{i=1}^n\theta(\widehat{X_i})\bigl(\rho(x_i)(a)\bigl)x_{n+1}\\
&~+\underline{a\sum_{i=1}^n\bigl[x_i,\mu(\widehat{X_i},x_{n+1})  \bigl]}+\rho(x_i)(a)\sum_{i=1}^n\mu\bigl(\widehat{X_i},x_{n+1}\bigl)+\sum_{i=1}^n\theta(\widehat{X_i})(a)[x_i,x_{n+1}]\\
&~+\sum_{i=1}^n\rho(x_i)\bigl(\theta(\widehat{X_i}(a))x_{n+1}\bigl)+\rho\circ\mu(X)(a)x_{n+1}+\underline{a\bigl[x_{n+1},\mu(X)\bigl]}\\
&=ad_{\ce}^n\mu(X,x_{n+1})+\alpha_{\mu,0}(\theta)(X)(a)x_{n+1},
\end{align*}\normalsize where the underlined terms corresponds to $ad_{\ce}^n\mu(X,x_{n+1})$ and
\begin{align*}
    \alpha_{\mu,0}(\theta)(X)(a)x_{n+1}&=\rho\circ\mu(X)(a)x_{n+1}+d_{\ce}^{n_1}\theta(X)(a)x_{n+1}\\
    &=\rho\circ\mu(X)(a)x_{n+1}+\sum_{1\leq i<j\leq n}\theta\bigl([x_i,x_j],\widehat{X_{ij}}\bigl)(a)x_{n+1}\\ &~+\sum_{i=1}^{n}\rho(x_i)\circ\theta(\widehat{X_i})(a)x_{n+1}+\sum_{i=1}^{n}\theta(\widehat{X_i})\circ\rho(x_i)(a)x_{n+1}. 
\end{align*} Thus, Eq. \eqref{eq:LRcohomop=2c} is satisfied. Next, let $x,z_2,\cdots,z_n\in L$ and $a\in A.$  We denote by $Z:=(z_2,\cdots,z_n)$, $\widehat{Z_i}:=(z_2,\cdots,\hat{z_i},\cdots,z_n)$ and $\widehat{Z_{ij}}:=(z_2,\cdots,\hat{z_i},\cdots,\hat{z_j},\cdots,z_n)$. We have
\begin{align*}
    &\,\,\,\delta^n(ax,Z)\\&=\bigl[ax,\mu(ax,Z)\bigl]+\sum_{i=2}^{n}\bigl[z_i,\w(az,\widehat{Z_i}) \bigl]+\mu\bigl((ax)^{[2]},Z \bigl)\\
    &~+\sum_{i=2}^{n}\mu\bigl([ax,z_i],ax,\widehat{Z_i} \bigl)+\sum_{1\leq i<j\leq n}\w\bigl( ax,[z_i,z_j],\widehat{Z_{ij}}\bigl)\\
    &=a\bigl[ax,\mu(x,Z)  \bigl]+\rho(ax)(a)\mu(x,Z)+\theta(Z)(a)[ax,x]+\rho(ax)\circ\theta(Z)(a)x\\
    &~+a^2\sum_{i=2}^{n}\bigl[z_i,\w(x,\widehat{Z_i}) \bigl]+\sum_{i=2}^{n}\theta(ax,\widehat{Z_i})(a)[z_i,x]+\sum_{i=2}^{n}\rho(z_i)\circ\theta(ax,\widehat{Z_i})(a)x\\
    &~+a^2\mu(x^{[2]},Z)+\rho(ax)(a)\mu(x,Z)+\theta(Z)\circ\rho(ax)(a)x\\
    &~+a\sum_{i=2}^{n}\mu\bigl([x,z_i],ax,\widehat{Z_i}\bigl)+\sum_{i=2}^{n}\theta(ax,\widehat{Z_i})(a)[x,z_i]+\sum_{i=2}^{n}\rho(z_i)(a)\mu(x,ax,\widehat{Z_i})\\
    &~+\sum_{i=2}^{n}\theta(ax,\widehat{Z_i})\circ\rho(z_i)(a)x+\sum_{2\leq i<j\leq n}a^2\w\bigl(x,[z_i,z_j],\widehat{Z_{ij}} \bigl)+\sum_{2\leq i<j\leq n}\theta\bigl(ax,[z_i,z_j],\widehat{Z_{ij}}\bigl)(a)x\\
    &=\underline{a^2\bigl[x,\mu(x,Z) \bigl]}+a\rho\circ\mu(x,Z)(a)x+\theta(Z)(a)\rho(x)(a)x+\underline{a\rho(x)\circ\theta(Z)(a)x}\\
    &~+\underline{a^2\sum_{i=2}^{n}\bigl[z_i,\w(z,\widehat{Z_i})\bigl]}+a\sum_{i=2}^{n}\rho(z_i)\circ\theta(x,\widehat{Z_i})(a)x+\underline{a^2\mu(x^{[2]},Z)}+\theta(Z)\circ\rho(ax)(a)x\\
    &~+\underline{a^2\sum_{i=2}^{n}\mu\bigl( [x,z_i],x,\widehat{Z_i}\bigl)}+\sum_{i=2}^{n}a\theta\bigl([x,z_i],\widehat{Z_i}\bigl)(a)x+\sum_{i=2}^{n}\rho(z_i)(a)\theta(x,\widehat{Z_i})(a)x\\
    &~+\sum_{i=2}^{n}\theta(x,\widehat{Z_i})\bigl(\rho(z_i)(a)\bigl)x+\underline{a^2\sum_{1\leq i<j\leq n}\w\bigl(x,[z_i,z_j],\widehat{Z_{ij}}\bigl)}+a\sum_{1\leq i<j\leq n}\theta\bigl(x,[z_i,z_j],\widehat{Z_{ij}}\bigl)(a)x\\
    &=a^2\delta^n(x,Z)+\alpha_{\mu,0}(\theta)(ax,Z)(a)x,
\end{align*}where the underlined terms corresponds to $a^2\delta^n(x,Z)$ and
\begin{align*}
    \alpha_{\mu,0}(\theta)(ax,Z)(a)x&=\rho\circ\mu(ax,Z)(a)x+\theta(Z)(a)\rho(x)(a)x+\sum_{1\leq i<j\leq n}\theta\bigl([z_i,z_j],ax,\widehat{Z_{ij}} \bigl)(a)x\\
    &~+\sum_{i=2}^{n}\theta\bigl([ax,z_i],\widehat{Z_{i}}\bigl)(a)x+\sum_{i=2}^{n}\bigl[\rho(z_i),\theta(ax,\widehat{Z_{i}})\bigl]_{\Der(A)}(a)x\\
    &~+\bigl[\rho(ax),\theta(Z)\bigl]_{\Der(A)}(a)x.
\end{align*}
Thus, Eq. \eqref{eq:LRcohomop=2d} is satisfied. We denote by $\widehat{Z_{ijk}}:=(z_2,\cdots,\hat{z_i},\cdots,\hat{z_j},\cdots,\hat{z_k},\cdots,z_n)$. We have 
\small{\begin{align*}
    &\delta\w(x,z_2,\cdots,az_i,\cdots,z_n)\\
    =&\bigl[x,a\mu(x,Z)\bigl]+\bigl[x,\theta(x,\widehat{Z_i})(a)z_i \bigl]+\sum_{\underset{j\neq i}{j=2}}^n\bigl[z_j,a\w(x,\widehat{Z_j})\bigl]+\sum_{\underset{j\neq i}{j=2}}^n\bigl[z_j,\gamma(x,\widehat{Z_{ij}})(a)z_i\bigl]\\
    &~+a\bigl[z_i,\w(x,\widehat{Z_i})\bigl]+\rho\circ\w(x,\widehat{Z_i})(a)z_i+a\mu(x^{[2]},Z)+\theta(x^{[2]},\widehat{Z_i})(a)z_i+\sum_{\underset{j\neq i}{j=2}}^na\mu\bigl([x,z_i],x,\widehat{Z_j}\bigl)\\  
    &~+\mu\bigl(a[x,z_i],x,\widehat{Z_i} \bigl)+\mu\bigl(\rho(x)(a)z_i,x,\widehat{Z_i} \bigl)+\sum_{\underset{j,k\neq i}{2\leq j<k\leq n}}a\w\bigl(x,[z_j,z_k],\widehat{Z_{jk}} \bigl)\\
    &~+\sum_{\underset{j,k\neq i}{2\leq j<k\leq n}}\gamma\bigl(x,[z_j,z_k],\widehat{Z_{ijk}}\bigl)(a)z_i+\sum_{\underset{j\neq i}{j=2}}^n\w\bigl(x,a[z_i,z_j],\widehat{Z_{ij}}\bigl)+\sum_{\underset{j\neq i}{j=2}}^n\w\bigl(x,\rho(z_j)(a)z_i,\widehat{Z_{ij}}\bigl)\\
    &=\underline{a\bigl[x,\mu(x,Z)\bigl]}+\rho(x)(a)\mu(x,Z)+\theta(x,\widehat{Z_i})(a)[x,z_i]+\rho(x)\bigl(\theta(x,\widehat{Z_i})(a) \bigl)\\
    &~+\underline{\sum_{\underset{j\neq i}{j=2}}^n a\bigl[z_j,\w(x,\widehat{Z_j}\bigl] }+\sum_{\underset{j\neq i}{j=2}}^n\rho(z_j)(a)\w(x,\widehat{Z_j})+\sum_{\underset{j\neq i}{j=2}}^n\gamma(x,\widehat{Z_{ij}})(a)[z_j,z_i]\\
    &~+\sum_{\underset{j\neq i}{j=2}}^n\rho(z_j)\bigl(\gamma(x,\widehat{Z_{ij}})(a)\bigl)(z_i)+\underline{a\bigl[z_i,\w(x,\widehat{Z_i})\bigl]}+\rho\circ\w(x,\widehat{Z_i})(a)z_i+\underline{a\mu(x^{[2]},Z)}+\theta(x^{[2]},\widehat{Z_i})(a)z_i\\
    &~+\underline{a\sum_{\underset{j\neq i}{j=2}}^n\mu\bigl([x,z_j],x,\widehat{Z_j}\bigl)}+\sum_{\underset{j\neq i}{j=2}}^n\theta\bigl([x,z_j],x,\widehat{Z_{ij}}\bigl)(a)z_i+\underline{a\mu\bigl([x,z_i],x,\widehat{Z_i}\bigl)}+\theta(x,\widehat{Z_i})(a)[x,z_i]\\
    &~+\rho(x)(a)\mu(z_i,x,\widehat{Z_i})+\theta(x,\widehat{Z_i}\bigl(\rho(x)(a)\bigl)z_i+\underline{a\sum_{\underset{j,k\neq i}{2\leq j<k\leq n}} \w\bigl(x,[z_j,z_k],\widehat{Z_{jk}}\bigl)  }\\
    &~+\sum_{\underset{j,k\neq i}{2\leq j<k\leq n}}\gamma\bigl(x,[z_j,z_k],\widehat{Z_{ijk}} \bigl)(a)z_i
    +\underline{a\sum_{\underset{j\neq i}{j=2}}^n\w\bigl(x,[z_i,z_j],\widehat{Z_{ij}} \bigl) }\\
    &~+\sum_{\underset{j\neq i}{j=2}}^n\gamma(x,\widehat{Z_{ij}})(a)[z_i,z_j]+\sum_{\underset{j\neq i}{j=2}}^n\rho(z_j)(a)\w(x,z_i,\widehat{Z_{ij}})+\sum_{\underset{j\neq i}{j=2}}^n\gamma(x,\widehat{Z_{ij}})\bigl(\rho(z_j)(a)\bigl)z_i\\
    &=a\delta^n\w(x,Z)+\beta_{\w,0}(\gamma)(x,\widehat{Z_i})(a)z_i,
\end{align*}}\normalsize{} where the underlined terms corresponds to $a\delta^n\w(x,Z)$ and 
\small{\begin{align*}
    &\beta_{\w,0}(\gamma)(x,\widehat{Z_i})(a)z_i\\=&~\rho\circ\gamma(x,\widehat{Z_i})(a)z_i+\bigl[\rho(x),\theta(x,\widehat{Z_i})\bigl]_{\Der(A)}(a)z_i+\sum_{\underset{j\neq i}{j=2}}^n\bigl[\rho(z_j),\gamma(x,\widehat{Z_{ij}})\bigl]_{\Der(A)}(a)z_i\\
    &~+\theta(x^{[2]},\widehat{Z_i})(a)z_i+\sum_{\underset{j\neq i}{j=2}}^n\theta\bigl([x,z_j],x,\widehat{Z_{ij}}\bigl)(a)z_i+\sum_{\underset{j,k\neq i}{2\leq j<k\leq n}}\gamma\bigl(x,[z_j,z_k],\widehat{Z_{ijk}}  \bigl)(a)z_i.
\end{align*}}\normalsize{}
Thus, Eq. \eqref{eq:LRcohomop=2e} is satisfied as well, which finishes the proof.
\end{proof}

We denote by 
\begin{align}
    Z^1_{\LR}(A,L,\rho)&:=\iota_1 \bigl(C^1_{\LR}(L;L)\bigl)\cap\mathfrak{Z}^1_{\mathrm{res}}(\rho,\rho);~B^1_{\LR}(A,L,\rho):=\iota_1 \bigl(C^1_{\LR}(L;L)\bigl)\cap\mathfrak{B}^1_{\mathrm{res}}(\rho,\rho);\\
    Z^2_{\LR}(A,L,\rho)&:=\iota_2 \bigl(C^2_{\LR}(L;L)\bigl)\cap\mathfrak{Z}^2_{\mathrm{res}}(\rho,\rho);~B^2_{\LR}(A,L,\rho):=\iota_2 \bigl(C^2_{\LR}(L;L)\bigl)\cap\mathfrak{B}^2_{\mathrm{res}}(\rho,\rho),
     \end{align}
and for $n\geq 3$,
     \begin{align}
       Z^n_{\LR}(A,L,\rho)&:=\iota_n \bigl(C^n_{\LR}(L;L)\bigl)\cap\mathfrak{Z}^n_{\mathrm{res}}(\rho,\rho);~B^n_{\LR}(A,L,\rho):=\iota_n \bigl(C^n_{\LR}(L;L)\bigl)\cap\mathfrak{B}^n_{\mathrm{res}}(\rho,\rho).     
\end{align}

\sssbegin{Theorem}\label{thm:p=2welldefined}
 Let $(A,L,\rho)$ be a restricted Lie Rinehart algebra. For all $n\geq 1$, the complex $\Bigl( C^n_{\LR}(L,L), \iota_{n+1}^{-1}\circ\fd^n_{\mathrm{res}}\circ\iota_n   \Bigl)$ is a well-defined cochain complex, called \textit{deformation cohomology complex} of $(A,L,\rho)$.
\end{Theorem}
\begin{proof}
    For all $n\geq 1$, Lemmas \ref{lem:diff1}, \ref{lem:diff2} and \ref{lem:diffn} ensure that the maps $\iota_{n+1}^{-1}\circ\fd^n_{\mathrm{res}}\circ\iota_n$ are well defined, that is, $\iota_{n+1}^{-1}\circ\fd^n_{\mathrm{res}}\circ\iota_n\bigl(C^{n}_{\LR}(A,L,\rho)\bigl)\subset C^{n+1}_{\LR}(A,L,\rho).$
    
    Moreover, they satisfy
    $$\iota_{n+2}^{-1}\circ\fd^{n+1}_{\mathrm{res}}\circ\iota_{n+1}\circ \iota_{n+1}^{-1}\circ\fd^n_{\mathrm{res}}\circ\iota_n=\iota_{n+2}^{-1}\circ\fd^{n+1}_{\mathrm{res}}\circ\fd^n_{\mathrm{res}}\circ\iota_n=0,$$ using Theorem \ref{thm:cohomorph2}. 
\end{proof}

\subsection{Formal deformations}\label{sec:deformations} Let $(A,L,\rho)$ be a restricted Lie-Rinehart algebra. Consider a triple $(\mu_t,\w_t,\rho_t)$ given by
\begin{equation*}
            \begin{array}{llllllll}
                \mu_t:&L\times L\longrightarrow L[[t]]& \text{ ;  }\quad & \w_t:&L\longrightarrow L[[t]]&\text{and}& \rho_t:&L\longrightarrow \Der\bigl(A[[t]]\bigl) \\[2mm]
                &(x,y)\longmapsto \displaystyle\sum_{i\geq 0}t^i \mu_i(x,y)&\quad  &&x\longmapsto \displaystyle\sum_{j\geq 0}t^j\omega_j(x)&&&x\longmapsto \displaystyle\sum_{k\geq 0}t^k\rho_k(x), \\[2mm]
             \end{array}
    \end{equation*}
where $(\mu_i,\w_i,\rho_i)\in C^2_{\LR}(L;L),~\forall i\geq 1$, $\mu_0=[\cdot,\cdot]$, $\w_0=(\cdot)^{[2]}$ and $\rho_0=\rho$. Such a triple is called \textit{formal deformation} of the restricted Lie-Rinehart algebra $(A,L,\rho)$ if the four following condition are satisfied, for all $x,y,z\in L:$
\begin{align}
    \underset{x,y,z}{\circlearrowleft}\mu_t\bigl(x,\mu_t(y,z)\bigl)&=0;\\
    \mu_t\bigl(x,\w_t(y)\bigl)&=\mu_t\bigl(\mu_t(x,y),y\bigl);\\
    \rho_t\bigl(\mu_t(x,y)\bigl)&=\rho_t(x)\circ\rho_t(y)+\rho_t(y)\circ\rho_t(x);\\
    \rho_t\bigl(\w_t(x)\bigl)&=\rho_t(x)\circ\rho_t(x).\end{align}
Note that since $(\mu_i,\w_i,\rho_i)\in C^2_{\LR}(L;L),~\forall i\geq 1$, the triple $(\mu_t,\w_t,\rho_t)$ automatically satisfies
    \begin{align}
        \mu_t(x,ay)&=a\mu_t(x,y)+\rho_t(x)(a)y,~\forall x,y\in L,~\forall a\in A; \text{ and }\\
        \w_t(ax)&=a^2\w_t(x)+\rho_t(ax)(a)x, ~\forall x\in L,~\forall a\in A.
    \end{align}
Applying \cite[Section 3.4.2]{EM2}, we have the following result.
\sssbegin{Proposition}
    Let $(\mu_t,\w_t,\rho_t)$ be a formal deformation of a restricted Lie-Rinehart algebra $(A,L,\rho)$. Then $(\mu_1,\w_1,\rho_1)\in  Z^2_{\LR}(A,L,\rho)$.
\end{Proposition}

Let $n\geq 1$. Consider a formal deformation $(\mu_t^n,\w_t^n,\rho_t^n)$ of order $n$ defined by
\begin{equation*}
            \begin{array}{llllllll}
                \mu_t^n:&L\times L\longrightarrow L[[t]]& \text{ ;  }\quad & \w_t^n:&L\longrightarrow L[[t]]&\text{and}& \rho_t^n:&L\longrightarrow \Der\bigl(A[[t]]\bigl) \\[2mm]
                &(x,y)\longmapsto \displaystyle\sum_{i= 0}^nt^i \mu_i(x,y)&\quad  &&x\longmapsto \displaystyle\sum_{j= 0}^nt^j\omega_j(x)&&&x\longmapsto \displaystyle\sum_{k=0}^nt^k\rho_k(x), \\[2mm]
             \end{array}
    \end{equation*}
We consider the problem of extending a deformation of order $n$ to a deformation of order $n+1$. Let
\begin{equation}\label{eq:ordern+1}
\mu_t^{n+1}:=\mu_t^n+t^{n+1}\mu_{n+1};~\w_t^{n+1}:=\w_t^n+t^{n+1}\w_{n+1};~\rho_t^{n+1}:=\rho_t^n+t^{n+1}\rho_{n+1},
\end{equation} where $\bigl(\mu_{n+1},\w_{n+1},\rho_{n+1}\bigl)\in C^2_{\LR}(L;L)$. Moreover, for $x,y,z\in L,$ we consider
\begin{align}
			\obs_{n+1}^{(1)}(x,y,z)&:=\underset{x,y,z}{\circlearrowleft}\sum_{i=1}^{n}\mu_i(x,\mu_{n+1-i}(y,z));\\ 
			\obs_{n+1}^{(2)}(x,y)&:=\sum_{i=1}^{n}\bigl( \mu_i(y,\w_{n+1-i}(x))+\mu_i(\mu_{n+1-i}(y,x),x) \bigl);\\
             \fobs^{(1)}_{n+1}(\rho)(x,y)&:=\sum_{i=1}^n\rho_i\circ\mu_{n+1-i}(x,y)+\sum_{\underset{i+j=n+1}{i,j\leq n}}\bigl[\rho_i(x),\rho_j(y)\bigl]_{\Der(A)};\\
    \fobs^{(2)}_{n+1}(\rho)(x)&:=\sum_{i=1}^n\rho_i\circ\w_{n+1-i}(x)+\sum_{\underset{0<i<j}{i+j=n+1}}\bigl[\rho_i(x),\rho_j(x)\bigl]_{\Der(A)}.
\end{align}
    
\sssbegin{Proposition}\label{prop:obs-p=2}
    Let $(\mu_t^n,\w_t^n,\rho_t^n)$ be a formal deformation of order $n$ of a restricted Lie-Rinehart algebra $(A,L,\rho)$. Moreover, suppose that the triple $(\mu_t^{n+1},\w_t^{n+1},\rho_t^{n+1})$ defined in Eq. \eqref{eq:ordern+1} is also a formal deformation of $(A,L,\rho)$. Then,
    \begin{equation} \bigl(\obs_{n+1}^{(1)},\obs_{n+1}^{(2)}\bigl)=d_{\mathrm{res}}^2\left(\mu_{n+1},~\w_{n+1} \right)
    \end{equation} and 
        \begin{equation} \bigl(\fobs_{n+1}^{(1)}(\rho),\fobs_{n+1}^{(2)}(\rho)\bigl)=\bigl(\alpha_{\mu_{n+1},0}(\rho_{n+1}),\beta_{\w_{n+1},0}(\rho_{n+1})\bigl)
        \end{equation} 
    \begin{proof}
        Using \cite[Lemma 4.20]{EM2} and \cite[Lemma 4.23]{EM2}, we have that
        $$ \Bigl(\bigl(\obs_{n+1}^{(1)},\obs_{n+1}^{(2)}\bigl),(0,0),\bigl(\fobs_{n+1}^{(1)}(\rho),\fobs_{n+1}^{(2)}(\rho)\bigl)\Bigl)\in \fC_{\mathrm{res}}^3(\rho,\rho). $$ Moreover, using \cite[Proposition 4.21]{EM2}, we have
        $$\bigl(\obs_{n+1}^{(1)},\obs_{n+1}^{(2)}\bigl)=d_{\mathrm{res}}^2\left(\mu_{n+1},~\w_{n+1} \right)$$ and using \cite[Proposition 4.24]{EM2}, we have
        $$\bigl(\fobs_{n+1}^{(1)}(\rho),\fobs_{n+1}^{(2)}(\rho)\bigl)=\bigl(\alpha_{\mu_{n+1},0}(\rho_{n+1}),\beta_{\w_{n+1},0}(\rho_{n+1})\bigl).$$
    \end{proof}
\end{Proposition}

%%%%%%%%%%%%%%%%%%%%%%%%%%%%%%%%%%%%%%%%%%%%%%
\subsection{Equivalence}\label{sec:equiv-2} Following \cite[Definition 3.14]{EM2}, two formal deformations $(\mu_t,\w_t,\rho_t)$ and $(\tilde{\mu_t},\tilde{\w_t},\tilde{\rho_t})$ of a restricted Lie-Rinehart algebra $(A,L,\rho)$ are \textit{equivalent} if there exists an $A$-linear formal automorphism $\phi_t:L[[t]]\rightarrow L[[t]]$ is defined on $L$ by 
	\begin{equation} 
    \phi_t(x)=\sum_{i\geq 0}t^i\phi_i(x),\text{ with }\phi_i:L\rightarrow L\text{ and }\phi_0=id,
        \end{equation}
    then extended by $\K[[t]]$-linearity, that satisfies for all $x,y\in L$
    \begin{align}
    \tilde{\mu_t}\bigl(\phi_t(x),\phi_t(y)\bigl)&=\phi_t\circ\mu_t(x,y);\\
        \tilde{\w_t}\circ\phi_t(x)&=\phi_t\circ\w_t(x),
    \end{align} and such that the following diagram commutes:
\begin{equation}\label{eq:diag-equiv-2}
 \begin{tikzcd}
{\bigl(L[[t]],\mu_t,\w_t\bigl)} \arrow[dd, "\phi_t"'] \arrow[rr, "\rho_t"]    &  & {\Der\bigl(A[[t]]\bigl)} \arrow[dd, "\id"] \\ &  &  \\
{\bigl(L[[t]],\tilde{\mu_t},\tilde{\w_t}\bigl)} \arrow[rr, "\widetilde{\rho_t}"'] &  & {\Der\bigl(A[[t]]\bigl)}.    
\end{tikzcd}
\end{equation}
Note that the diagram \eqref{eq:diag-equiv-2} is equivalent to  $\tilde{\rho_t}\circ\phi_t=\rho_t.$ Comparing with Equation \eqref{eq:morphRLR}, we have that $\phi_t$ is a morphism of restricted Lie-Rinehart algebras. 
\sssbegin{Proposition}\label{prop:equiv-p=2}
    Let $(\mu_t,\w_t,\rho_t)$ and $(\tilde{\mu_t},\tilde{\w_t},\tilde{\rho_t})$ be two equivalent deformations of a restricted Lie-Rinehart algebra $(A,L,\rho)$. Then, the infinitesimals $(\mu_1,\w_1,\rho_1)$ and $(\tilde{\mu_1},\tilde{\w_1},\tilde{\rho_1})$ differ from an element of $B^2_{\LR}(A,L,\rho)$.
\end{Proposition}
\begin{proof}
Using \cite[Proposition 4.16]{EM2}, we have that
    $\bigl(\mu_1-\tilde{\mu_1},\w_1-\tilde{\w_1}\bigl)\in B_{\mathrm{res}}^2(L,L).$
    Moreover, \eqref{eq:diag-equiv-2} leads to 
       $\rho_1-\tilde{\rho_1}=\alpha_{\phi_1,0}(0).$
    Thus, we have \begin{equation}\bigl(\mu_1-\tilde{\mu_1},\w_1-\tilde{\w_1},\rho_1-\tilde{\rho_1}\bigl)\in \mathfrak{B}^2_{\mathrm{res}}(\rho,\rho)\end{equation} and the result follows.
\end{proof}
\subsubsection{Trivial deformations} A deformation $(\mu_t,\w_t,\rho_t)$ of a restricted Lie-Rinehart algebra $(A,L,\rho)$ is called \textit{trivial} if it is equivalent to the deformation given by $\tilde{\mu_t}=[\cdot,\cdot],~\tilde{\w_t}=(\cdot)^{[2]},~\tilde{\rho_t}=\rho.$ Namely, we have
\begin{equation}
    \phi_t\circ\mu_t(x,y)=[\phi_t(x),\phi_t(y)];~\phi_t\circ\w_t(x)=\phi_t(x)^{[2]};~\rho_t(x)=\rho\circ\phi_t(x),~\forall x\in L.
\end{equation} Then, we have $(\mu_1,\w_1,\rho_1)\in B^2_{\LR}(A,L,\rho)$. Therefore, any non trivial deformation is equivalent to a deformation such that there exists $n\in\mathbb{N}$,
$$(\mu_i,\w_i,\rho_i)\in B^2_{\LR}(A,L,\rho)~\forall i<n, \text{ and } (\mu_n,\w_n,\rho_n)\in Z^2_{\LR}(A,L,\rho)\setminus B^2_{\LR}(A,L,\rho). $$
\subsection{Examples}\label{sec:ex}
We provide some examples of computation of restricted Lie-Rinehart deformation 2-cocycles. 

\subsubsection{Associative commutative algebras of dimension 2 in characteristic 2} The following table displays the classification of associative commutative algebras of dimension 2 in characteristic 2, see  \cite[Theorem 1]{GK}, as well as their derivations. The non-written products are zero.

\begin{center}
	\begin{tabular}{|c|c|c|}
		\hline
		Name&Associative commutative product& Derivations \\\hline
		$A_1$&$e_1e_1=e_1$& $e_2\otimes e_2^*$  \\\hline
		$A_2$&$e_1e_1=e_2$& $e_2\otimes e_1^*$ \\\hline
		$A_3$&$e_1e_1=e_1;~e_2e_2=e_2$&$0$ \\\hline
		$A_4$&$e_1e_1=e_1;~e_1e_2=e_2;~e_2e_1=e_2$& $e_1\otimes e_2^*;e_2\otimes e_2^*$  \\\hline	$A_5$&$e_1e_1=e_1;~e_1e_2=e_2;~e_2e_1=e_2;~e_2e_2=e_1+e_2$& $0$\\\hline
	\end{tabular}
	\\~\\ Associative commutative algebras of dimension 2.
\end{center}

\subsubsection{A rigid restricted Lie-Rinehart algebra} Let $L$ be the 2 dimensional restricted Lie algebra spanned by elements $x,y$ with the bracket $[x,y]=y$ and the 2-map $x^{[2]}=x,~y^{[2]}=0$. The associative commutative algebra $A_4$ acts on $L$ by
\begin{equation}
    e_1\cdot x=x,~e_1\cdot y=y,~e_2\cdot x=y,~e_2\cdot y=0.    
\end{equation}
Consider the $A$-linear restricted Lie algebras morphism $\rho:L\rightarrow\Der(A_4)$ given by
\begin{equation}
    \rho(x)(e_1)=\rho(y)(e_1)=0,~\rho(x)(e_2)=e_2,~\rho(y)(e_2)=0. \end{equation}
Then, $\rho$ satisfies Eqs. \eqref{eq:leib}-\eqref{eq:hoch} and $(A_4,L,\rho)$ is a restricted Lie-Rinehart algebra. We have $Z^{2}_{\mathrm{res}}(L;L)=\{0\}$. Let $\theta:L\rightarrow \Der(A)$. Eq. \eqref{eq:LRcohomop=2a} leads to $\theta=0$. Therefore, $Z_{\LR}^2(A,L,\rho)=\{0\}.$

\subsubsection{The 2-dimensional abelian Lie algebra} Let $L_{\ab}$ be the two dimensional abelian Lie algebra spanned by elements $x,y$. The associative commutative algebra $A_4$ acts on $L_{\ab}$ by
\begin{equation}
    e_1\cdot x=x,~e_1\cdot y=y,~e_2\cdot x=y,~e_2\cdot y=0.    
\end{equation} We consider the 2-maps $(\cdot)^{[2]_0}= 0$ and $x^{[2]_1}=\lambda_1 x+\lambda_2 y,~\lambda_1,\lambda_2\in\K, y^{[2]_1}=0$. We denote the resulting restricted Lie algebras by $L_{\ab}^0$ and $L_{\ab}^1$, respectively. It is not difficult to see that any anchor map $\rho:L_{\ab}^i\rightarrow\Der(A),~i\in\{0,1\},$ must vanish. Therefore, we consider the restricted Lie-Rinehart algebras $(A_4,L_{\ab}^0,0)$ and $(A_4,L_{\ab}^1,0)$.

\paragraph{The restricted Lie-Rinehart algebra $(A_4,L_{\ab}^0,0)$} A basis of $Z_{\mathrm{res}}^2(L_{\ab}^0;L_{\ab}^0)$ is given by
\begin{equation}
    \Bigl\{(\mu_1,0),~(\mu_1,0),~(0,\w_1),~(0,\w_2),~(0,\w_3),~(0,\w_4)   \Bigl\},
\end{equation} where $\mu_1=x\otimes x^*\wedge y^*,~\mu_2=y\otimes x^*\wedge y^*$ and $\w_1(x)=x,~\w_2(x)=y,~\w_3(y)=x,~\w_4(y)=y,$ the other images being zero. Consider an $A$-linear map $\theta_1:L_{\ab}^0\rightarrow\Der(A)$ Then, $(\mu_1,\theta_1)$ satisfy Eq. \eqref{eq:LRcohomop=2a} if and only if
\begin{equation}
    \theta_1(x)(e_1)=\theta_1(y)(e_1)=0;~\theta_1(x)(e_2)=e_1;~\theta_1(y)(e_2)=e_2.
\end{equation} Moreover, $(\w_i,\mu_1)$ satisfy Eq. \eqref{eq:LRcohomop=2b} if and only if $i=4$. Finally, we have $\alpha_{\mu_1,0}(\theta_1)=0$ and $\beta_{\w_4,0}(\theta_1)=0.$ Thus, $(\mu_1,\w_4,\theta_1)\in Z_{\LR}^2(A_4,L_{\ab}^0,0)$. Similarly, we have that $$(\mu_2,\w_1,\theta_2)\in Z_{\LR}^2(A_4,L_{\ab}^0,0) \text{ and } (\mu_2,\w_2,\theta_2)\in Z_{\LR}^2(A_4,L_{\ab}^0,0),$$ with $\theta_2:L_{\ab}^0\rightarrow\Der(A)$ given by
\begin{equation}
    \theta_2(x)(e_1)=\theta_2(y)(e_1)=0;~\theta_2(x)(e_2)=e_2;~\theta_2(y)(e_2)=0.
\end{equation} Since $(\w_1,0)$ and $(\w_2,0)$ also satisfy Eq. \eqref{eq:LRcohomop=2b}, a basis of the space $ Z_{\LR}^2(A_4,L_{\ab}^0,0)$ is given by
\begin{equation}
 \Bigl\{(\mu_1,\w_4,\theta_1),~(\mu_2,0,\theta_2),~(0,\w_1,0),~(0,\w_2,0)  \Bigl\}.
\end{equation} None of those cocycles are coboundaries, since we have $\fd^1_{\mathrm{res}}=0$ in that case.

\paragraph{The restricted Lie-Rinehart algebra $(A_4,L_{\ab}^1,0)$} A basis of $Z_{\mathrm{res}}^2(L_{\ab}^1;L_{\ab}^1)$ is given by
\begin{equation}
    \Bigl\{(\mu_1+\mu_2,0),~(0,\w_1),~(0,\w_2),~(0,\w_3),~(0,\w_4)   \Bigl\},
\end{equation} with the same notations as in the previous example. Similarly, a basis of the space $ Z_{\LR}^2(A_4,L_{\ab}^1,0)$ is given by
\begin{equation}
 \Bigl\{(\mu_1+\mu_2,\w_4,\theta_1),~(0,\w_1,0),~(0,\w_2,0)  \Bigl\}.
\end{equation}None of those cocycles are coboundaries.\\

%%%%%%%%%%%%%%%%%%%%%%%%%%%%%%%%%%%%%%%%%%%%%%%%%%%%%%%%%%%%%%%%%%%%%%%%%%%%%%%
\newpage
\section{Deformations of restricted Lie-Rinehart algebras, $p\geq 3$}\label{sec:p-geq3}

In this section, $\K$ is a field of characteristic $p\geq 3$. In that case, the situation is much more complicated due to the presence of the $p-1$ power in Eq. \eqref{eq:hoch}. We construct deformation 2-cocycles and explore their connections to formal deformations.

\subsection{Cohomological constructions} We introduce the cohomological objects involved in the deformation theory of restricted Lie-Rinehart algebras.
\subsubsection{Restricted cohomology of restricted morphisms} We recall the construction of the restricted cohomology of restricted morphisms developed in \cite{EM2}. For the definition of the restricted cohomology for restricted Lie algebras introduced by Evans-Fuchs, see \cite{EM2}.\footnote{ Those formulas were introduced in \cite{EF}, but we refer to \cite{EM2} for consistency with the present paper. Moreover, the spaces $C^k_{\mathrm{res}}$ are denoted by $C^k_{\mathrm{*}}$ in \cite{EM2}.} Let $(L,[ \cdot , \cdot ]_L,(\cdot )^{[p]_L})$ and $(M,[ \cdot , \cdot ]_M,(\cdot )^{[p]_M})$ be restricted Lie algebras and let $\varphi:L\rightarrow M$ be a restricted morphism. The restricted Lie algebra $M$ has a restricted $L$-module structure given by $x\cdot m:=[\varphi(x),m]_M$.
The restricted cochain spaces are defined by
	\begin{align*}
		\fC^0_{\mathrm{res}}(\varphi,\varphi)&:=0;\\
		\fC^1_{\mathrm{res}}(\varphi,\varphi)&:=C^1_{\mathrm{res}}(L;L)\times  C^1_{\mathrm{res}}(M;M)\times C^0_{\mathrm{res}}(L;M);\\
		\fC^2_{\mathrm{res}}(\varphi,\varphi)&:=C^2_{\mathrm{res}}(L;L)\times  C^2_{\mathrm{res}}(M;M)\times C^1_{\mathrm{res}}(L;M);\\
		\fC^3_{\mathrm{res}}(\varphi,\varphi)&:=C^3_{\mathrm{res}}(L;L)\times  C^3_{\mathrm{res}}(M;M)\times \widetilde{C}^2_{\mathrm{res}}(L;M),
	\end{align*}
	where $\widetilde{C}^2_{\mathrm{res}}(L;M):=\bigl\lbrace (\psi,\w),~\psi\in C^2_{\mathrm{CE}}(L,M),~\w: L\rightarrow M,~\w \text{ is $p$-homogeneous}\bigl\rbrace.$\footnote{Recall that such a map $\w$ is called $p$-homogeneous if $\w(\lambda x)=\lambda^p\w(x),\quad\forall \lambda\in\K,~\forall x\in L$.} We define the restricted differentials 		
	\begin{align*}		\fd_{\mathrm{res}}^1:~&\fC_{\mathrm{res}}^1(L;M)\rightarrow \fC^2_{\mathrm{res}}(L;M),\\\fd_{\mathrm{res}}^1(\gamma,\tau,m)&:=\Bigl(\bigl(d_{\mathrm{CE}}^1\gamma,\ind^1(\gamma)\bigl),\bigl(d_{\mathrm{CE}}^1\tau,\ind^1(\tau)\bigl),\varphi\circ\gamma-\tau\circ\varphi-[\varphi(\cdot),m]_M\Bigl);\\
		\fd_{\mathrm{res}}^2:~&\fC_{\mathrm{res}}^2(L;M)\rightarrow \fC^3_{\mathrm{res}}(L;M),\\\fd_{\mathrm{res}}^2\bigl((\mu,\omega),(\nu,\epsilon),\theta\bigl)&:=\Bigl( \bigl(d_{\mathrm{CE}}^2\mu,\ind^2(\mu,\w)\bigl),\bigl(d_{\mathrm{CE}}^2\nu,\ind^2(\nu,\epsilon)\bigl),\bigl(\alpha_{\mu,\nu}(\theta),\beta_{\w,\epsilon}(\theta)\bigl)\Bigl),
	\end{align*}
	where \begin{equation}\alpha_{\mu,\nu}(\theta):=\varphi\circ\mu-\nu\circ(\varphi^{\otimes 2})-d^{1}_{\ce}\theta,\end{equation} 
    and \begin{equation}\beta_{\w,\epsilon}(\theta)(x):=\theta(x^{[p]})+\varphi(\w(x))-\epsilon(\varphi(x))-x^{p-1}\cdot\theta(x),~\forall x\in L.\end{equation} We denote $\mathfrak{Z}^n_{\mathrm{res}}(\varphi,\varphi):=\Ker(\fd^n_{\mathrm{res}})$ and $\mathfrak{B}^n_{\mathrm{res}}(\varphi,\varphi):=\text{Im}(\fd^{n-1}_{\mathrm{res}})$, for $n=1,2$.

\subsubsection{Deformation 2-cocycles}
Let $(A,L,\rho)$ be a restricted Lie-Rinehart algebra (see Definition \ref{def:RLRp}). 
For every derivation $d\in\Der(A),$ we define
\begin{equation}
    C_{\LR}^1(L,L)^{d}:=\Bigl\{(\gamma,d),~ \gamma\in C_{\mathrm{res}}^1(L;L), \gamma(ax)=a^p\gamma(x)+d(a)x,~\forall x\in Lm~\forall a\in A\Bigl\}.
    \end{equation}
    There is an embedding
    \begin{equation}
        \iota_1:C^1_{\LR}(L;L)^{d}~\hookrightarrow\fC^1_{\mathrm{res}}(\rho,\rho),\quad (\gamma,d)\mapsto (\gamma,0,d).
    \end{equation}
For every $A$-linear map $\theta:L\rightarrow\Der(A)$, we define
\begin{equation}
    C_{\LR}^2(L,L)^{\theta}:=\Bigl\{(\mu,\w,\theta),~(\mu,\w)\in C_{\mathrm{res}}^2(L;L) \text{ satisfying } \eqref{eq:LRcohomop=2a} \Bigl\}.
    \end{equation}
    There is an embedding
    \begin{equation}
        \iota_2:C^2_{\LR}(L;L)^{\theta}~\hookrightarrow\fC^2_{\mathrm{res}}(\rho,\rho),\quad (\mu,\w,\theta)\mapsto \Bigl((\mu,\w),(0,0),\theta \Bigl).
    \end{equation}
    Finally, let
\begin{equation}
   C^1_{\LR}(L;L):=\bigcup_{d}C^1_{\LR}(L;L)^{d} ~\text{ and }~C^2_{\LR}(L;L):=\bigcup_{\theta}C^2_{\LR}(L;L)^{\theta}.
\end{equation} We are now in position to define
\begin{equation}
  Z^2_{\LR}(A,L,\rho):=\Bigl\{ (\mu,\w,\theta)\in  C^2_{\LR}(A,L,\rho) \text{ satisfying } \fd^2_{\mathrm{res}}\circ\iota_2(\mu,\w,\theta)=0 \text{ and } \eqref{eq:cohomopp}   \Bigl\},
\end{equation} where
\begin{equation}\label{eq:cohomopp}
\w(ax)-a^p\w(x)=a^{p-1}\sum_{i=0}^{p-2}\rho(x)^i\circ\theta(x)\circ\rho(x)^{p-2-i}(a)x,~\forall x\in L,~\forall a\in A.
\end{equation}
Elements of  $Z^2_{\LR}(A,L,\rho)$ are called \textit{deformation 2-cocycles} of $(A,L,\rho)$. A deformation cocycle $(\mu,\w,\theta)\in Z^2_{\LR}(A,L,\rho)$ is called \textit{trivial} if there exists $(\gamma,d)\in C^1_{\LR}(L;L)$ such that $\iota_2(\mu,\w,\theta)=\fd_{\mathrm{res}}^1\iota_1(\gamma,d)$.
\subsection{Formal deformations} Consider a triple $(\mu_t,\w_t,\rho_t)$ given by
\begin{equation*}
            \begin{array}{llllllll}
                \mu_t:&L\times L\longrightarrow L[[t]]& \text{ ;  }\quad & \w_t:&L\longrightarrow L[[t]]&\text{and}& \rho_t:&L\longrightarrow \Der\bigl(A[[t]]\bigl) \\[2mm]
                &(x,y)\longmapsto \displaystyle\sum_{i\geq 0}t^i \mu_i(x,y)&\quad  &&x\longmapsto \displaystyle\sum_{j\geq 0}t^j\omega_j(x)&&&x\longmapsto \displaystyle\sum_{k\geq 0}t^k\rho_k(x), \\[2mm]
             \end{array}
    \end{equation*}
where $(\mu_i,\w_i,\rho_i)\in C^2_{\LR}(L;L),~\forall i\geq 1$, $\mu_0=[\cdot,\cdot]$, $\w_0=(\cdot)^{[p]}$ and $\rho_0=\rho$. Such a triple is called \textit{formal deformation} of the restricted Lie-Rinehart algebra $(A,L,\rho)$ if the five following condition are satisfied, for all $x,y,z\in L$ and for all $a\in A:$
\begin{align}
    \underset{x,y,z}{\circlearrowleft}\mu_t\bigl(x,\mu_t(y,z)\bigl)&=0;\\
\mu_t\bigl(x,\w_t(y)\bigl)&=\mu_t\bigl(\mu_t(\cdots,\mu_t(x,\underset{p \text{ terms}}{\underbrace{y),\cdots y),y}}\bigl);\\\label{eq:hoch-t}
    \rho_t\bigl(\mu_t(x,y)\bigl)&=\rho_t(x)\circ\rho_t(y)+\rho_t(y)\circ\rho_t(x);\\
    \rho_t\bigl(\w_t(x)\bigl)&=\rho_t(x)^p;\\\label{eq:fifth}
    \w_t(ax)&=a^p\w_t(x)+\rho_t(ax)^{p-1}(a)x.
\end{align}
Note that since $(\mu_i,\w_i,\rho_i)\in C^2_{\LR}(L;L),~\forall i\geq 1$, we automatically have
\begin{equation}
\mu_t(x,ay)=a\mu_t(x,y)+\rho_t(x)(a)y,~\forall x,y\in L,~\forall a\in A.
\end{equation}
In contrast to the $ p = 2 $ case, a fifth deformation equation (Eq.~\eqref{eq:fifth}) arises due to the presence of the $ p-1 $ power in Eq.~\eqref{eq:hoch}, which leads to nontrivial relations between the components of $ \rho_t $.

\sssbegin{Proposition}\label{prop:cocy-p}
    Let $(\mu_t,\w_t,\rho_t)$ be a deformation of a restricted Lie-Rinehart algebra $(A,L,\rho)$. Then, we have $(\mu_1,\w_1,\rho_1)\in Z^2_{\LR}(A,L,\rho)$.
\end{Proposition}

\begin{proof}
    By \cite[Section 3.4.2]{EM2}, we have $\iota_2(\mu_1,\w_1,\rho_1)\in \fZ^2_{\mathrm{res}}(\rho,\rho)$. Moreover, expanding Eq. \eqref{eq:hoch-t} for $x\in L$ and $a\in A$ gives
    \begin{equation}
        \sum_{k\geq 0}t^k\w_k(ax)=\sum_{k\geq 0}t^k a^p\w(x)+a^{p-1} \sum_{k\geq 0}t^k \sum_{i_1+\cdots+i_{p-1}=k} \rho_{i_1}(x)\circ\cdots\circ\rho_{i_{p-1}}(x)(a)y
    \end{equation}
    Collecting the coefficients of $t$ gives 
    \begin{equation}
       \w_1(ax)=a^p\w_1(x)+a^{p-1}\sum_{i=0}^{p-2}\rho(x)^i\circ\rho_1(x)\circ\rho(x)^{p-2-i}(a)x,~\forall x\in L,~\forall a\in A,
    \end{equation} which finishes the proof.
\end{proof}
Next, we investigate the problem of extending a deformation of order 1 to a deformation of order 2. Let 
\begin{equation}\label{eq:order2}
\mu_t^{2}:=\mu_t^1+t^{2}\mu_{2};~\w_t^{2}:=\w_t^1+t^{2}\w_{2};~\rho_t^{2}:=\rho_t^1+t^{2}\rho_{2},
\end{equation} where the triple $\bigl(\mu_{2},\w_{2},\rho_{2}\bigl)\in C^2_{\LR}(L;L)$ and $\bigl(\mu_{t}^1,\w_{t}^1,\rho_{t}^1\bigl)$ is a formal deformation of order 1 of $(A,L,\rho)$.

For $x,y,z\in L,$ we consider
\begin{align}
			\obs_{2}^{(1)}(x,y,z)&:=-\underset{x,y,z}{\circlearrowleft}\mu_1(x,\mu_{1}(y,z));\\ 
			\obs_{2}^{(2)}(x,y)&:= \mu_1(x,\w_{1}(y))-\mu_1(\mu_{1}(x,y),y) ;\\
             \fobs^{(1)}_{2}(\rho)(x,y)&:=-\rho\bigl(\mu_1(x,y)\bigl)+\bigl[\rho_1(x),\rho_1(y)\bigl]_{\Der(A)};\\    
    \fobs^{(2)}_{2}(\rho)(x)&:=\ad_{\rho(x)}^{p-1}\circ\rho_1(x)-\rho\circ\w_1(x)\\\nonumber&~~-\sum_{i+j=p-2}\ad^i_{\rho(x)}\Bigl(\bigl[\rho_1(x),\ad^j_{\rho(x)}\circ\rho_1(x)  \bigl]_{\Der(A)}\Bigl).
\end{align}

\sssbegin{Proposition}\label{prop:obs-p}
Let $\bigl(\mu_{t}^1,\w_{t}^1,\rho_{t}^1\bigl)$ be a formal deformation of order 1 of a restricted Lie-Rinehart algebra $(A,L,\rho)$. Suppose that $\bigl(\mu_{t}^2,\w_{t}^2,\rho_{t}^2\bigl)$ defined in Eq. \eqref{eq:order2} is also a formal deformation of $(A,L,\rho)$. Then,
\begin{align}\label{eq:obs.restrp}
    \bigl(\obs_2^{(1)},\obs_2^{(2)}\bigl)&=d_{\mathrm{res}}^2(\mu_2,\w_2),\\
\label{eq:obs.morphp}
    \fobs_2^{(1)}(\rho)&=-\alpha_{\mu_1,0}(\rho_2);~ \fobs_2^{(2)}(\rho)=-\beta_{\w_1,0}(\rho_2), 
\end{align} and moreover
\begin{align}\label{eq:eqenplusobsp}
    \w_2(ax)-a^p\w_2(x)=~&a^{p-1}\sum_{i=0}^{p-2}\rho(x)\circ\rho_2(x)\circ\rho(x)^{p-2-i}(a)x\\\nonumber
    &~+a^{p-1}\sum_{i+j+k=p-3}\rho(x)^i\circ\rho_1(x)\circ\rho(x)^j\circ\rho_1(x)\circ\rho(x)^k(a)x.
\end{align}
\end{Proposition}

\begin{proof}
    Equations \eqref{eq:obs.restrp} and \eqref{eq:obs.morphp} are consequences of \cite[Proposition 3.10]{EM2} and \cite[Proposition 3.13]{EM2}, respectively. Moreover, for $x\in L$ and $a\in A$, expanding $$\w_t^2(ax)=a^p\w_t^2(x)+\rho_t^2(ax)^{p-1}(a)x$$ and collecting the coefficients of $t^2$ leads to Eq. \eqref{eq:eqenplusobsp}.
\end{proof}

 \noindent\textbf{Equivalence.} Equivalences of formal deformations of restricted Lie-Rinehart algebras in characteristic $p\geq3$ can be investigated similarly to the $p=2$ case, see Section \ref{sec:equiv-2}. In particular, we have the following.

 \sssbegin{Proposition}\label{prop:equiv-pgeq3}
    Let $(\mu_t,\w_t,\rho_t)$ and $(\tilde{\mu_t},\tilde{\w_t},\tilde{\rho_t})$ be two equivalent deformations of a restricted Lie-Rinehart algebra $(A,L,\rho)$. Then, $(\mu_1-\tilde{\mu_1},\w_1-\tilde{\w_1},\rho_1-\tilde{\rho_1})$ is a trivial deformation cocycle.
\end{Proposition}

\noindent\textbf{Acknowledgments.} The author thanks the referee for his suggestions, as well as S. Bouarroudj and A. Makhlouf for their support, encouragement and many fruitful discussions.
%%%%%%%%%%%%%%%%%%%%%%%%%%%%%%%%%%%%%%%%%%%%%%%%%%%%%%%%%%%%%%%%%%%%%%%%%%%%%%

\end{document}